\documentclass[11pt]{article}

\usepackage{fullwidth}
\usepackage{color}
\usepackage[color,matrix,arrow]{xy}
\usepackage[top=1.5in, bottom=1.5in, left=1in, right=1in]{geometry}
\usepackage{amsgen}
\usepackage{amsmath}
\usepackage{amstext}
\usepackage{amsbsy}
\usepackage{amsopn}
\usepackage{amsfonts}
\usepackage{amssymb}
\usepackage{eepic}
\usepackage{graphicx}
\usepackage{tikz}
\usepackage{epsf}
\usepackage{pstricks}
\usepackage{mathrsfs}
\usepackage{faktor}
\usepackage{amsthm}
\xyoption{all}

\def\Box{\square}

\def\tra#1{\smash{\mathop{\mid\kern
-1pt\joinrel\relbar\joinrel\relbar}\limits^{*}_{#1}}}
\def\longtra#1{\smash{\mathop{\mid\kern
-1pt\joinrel\relbar\joinrel\relbar\joinrel\relbar}\limits^{*}_{#1}}}
\def\vlongtra#1{\smash{\mathop{\mid\kern
-1pt\joinrel\relbar\joinrel\relbar\joinrel\relbar\joinrel\relbar}\limits^{*}_{#1}}}
\def\vvlongtra#1{\smash{\mathop{\mid\kern
-1pt\joinrel\relbar\joinrel\relbar\joinrel\relbar\joinrel\relbar\joinrel\relbar}\limits^{*}_{#1}}}
\def\vvvlongtra#1{\smash{\mathop{\mid\kern
-1pt\joinrel\relbar\joinrel\relbar\joinrel\relbar\joinrel\relbar\joinrel\relbar\joinrel\relbar}\limits^{*}_{#1}}}
\def\etra#1{\smash{\mathop{\mid\kern
-1pt\joinrel\relbar\joinrel\relbar}\limits_{#1}}}

\def\max{\mbox{max}}

\def\bi{\begin{itemize}}
\def\ei{\end{itemize}}
\def\beq{\begin{equation}}
\def\eeq{\end{equation}}

\def\J{{\cal{J}}}

\def\GM{\operatorname{GM}}
\def\GGM{\operatorname{GGM}}

\def\RLM{\operatorname{RLM}}
\def\RM{\operatorname{RM}}
\def\LM{\operatorname{LM}}

\theoremstyle{plain}
\newtheorem{T}{Theorem}[section]
\newcommand{\bt}{\begin{T}}
\newcommand{\et}{\end{T}}
\newcommand{\ftd}{$\square$\end{T}}

\newtheorem{Proposition}[T]{Proposition}
\newcommand{\bp}{\begin{Proposition}}
\newcommand{\ep}{\end{Proposition}}
\newcommand{\fpd}{$\square$\end{Proposition}}

\newtheorem{Lemma}[T]{Lemma}
\newcommand{\bl}{\begin{Lemma}}
\newcommand{\el}{\end{Lemma}}
\newcommand{\fld}{$\square$\end{Lemma}}

\newtheorem{Corol}[T]{Corollary}
\newcommand{\bc}{\begin{Corol}}
\newcommand{\ec}{\end{Corol}}
\newcommand{\fcd}{$\square$\end{Corol}}

\newtheorem{Result}[T]{Result}
\newcommand{\br}{\begin{Result}}
\newcommand{\er}{\end{Result}}
\newcommand{\frd}{$\square$\end{Result}}

\theoremstyle{definition}
\newtheorem{Example}[T]{Example}
\newcommand{\be}{\begin{Example}}
\newcommand{\ee}{\end{Example}}

\newtheorem{Problem}[T]{Problem}
\newcommand{\bq}{\begin{Problem}}
\newcommand{\eq}{\end{Problem}}

\newtheorem{Remark}[T]{Remark}
\newcommand{\brm}{\begin{Remark}}
\newcommand{\erm}{\end{Remark}}

\newtheorem{Definition}[T]{Definition}
\newcommand{\bd}{\begin{Definition}}
\newcommand{\ed}{\end{Definition}}

\newtheorem{Construction}[T]{Construction}
\newcommand{\bco}{\begin{Construction}}
\newcommand{\eco}{\end{Construction}}

\newtheorem{Computation}[T]{Computation}
\newcommand{\bcomp}{\begin{Computation}}
\newcommand{\ecomp}{\end{Computation}}

\newcommand{\Qed}{\hfill
$\Box$
\par\bigbreak}

\normalbaselines
\def\abstract#1{\par\bigskip
\begingroup\small
\baselineskip=12truept
\begin{center}ABSTRACT\end{center}
\par\medskip\par\noindent
\null\hfill\hbox{\vbox{\hsize=5truein\noindent#1}}
\hfill\null\par\endgroup\par}

\title{Aperiodic Flows on Finite Semigroups: Foundations and First Examples}
\author{Stuart Margolis and John Rhodes \\ \\
Dedicated to the Memory of Professor Boris I. Plotkin}

\date{}

\begin{document}

\maketitle

\abstract{
The theory of flows is used as a crucial tool in the recent proof by Margolis, Rhodes and Schilling that Krohn-Rhodes complexity is decidable. In this paper we begin a systematic study of aperiodic flows. We give the foundations of the theory of flows and give a unified approach to the Presentation Lemma and its relations to flows and the Slice Theorem. We completely describe semigroups having a flow over the trivial semigroup and connect this to classical results in inverse semigroup theory. We reinterpret Tilson's Theorem on the complexity of small monoids in terms of flows. We conclude with examples of semigroups built from the character table of Abelian Groups that have an aperiodic flows.}


\section{Introduction}\label{Intro}
In the papers \cite{complexity1} and \cite{complexityn}  Margolis, Rhodes and  Schilling proved that it is decidable to compute the complexity of a finite semigroup. This solved one of the major problems of finite semigroup theory. A key tool in the proof is the notion of a flow on an automaton. This was first defined and developed in \cite{Trans}. In particular a $\GM$ semigroup has complexity 1 if and only if its $\RLM$ image has complexity at most 1 and there is a flow from an aperiodic transformation semigroup to the Rhodes lattice associated to $S$. The purpose of this paper is to begin an in depth study of aperiodic flows as objects of interest by themselves. We illustrate flow theory by computing examples of flows for some classes of semigroups. This will be the first paper of a series of papers on flow theory \cite{FlowsII, FlowsIII}. Here is an outline of the paper.

Sections 2 and 3 give background material that we include for the benefit of the reader. Section 2 is devoted to semilocal theory  the type II subsemigroup and the Tilson congruence on a transformation semigroup. More details can be found in Chapters 7-8 of \cite{Arbib}, Chapter 4 of \cite{qtheory} and \cite{Redux}. Section 3 describes the connection between two classes of lattices that arise in flow theory. The first is the set-partition lattice $SP(X)$ on a set $X$ and the second are the Rhodes lattices $Rh_{B}(G)$. 
Section 4 gives the definition of flow from \cite{Trans} and summarizes the connections of flows to the Presentation Lemma \cite{AHNR.1995} and \cite[4.14]{qtheory}. One of the important goals of this paper is to give a unified version of the Presentation Lemma and flow theory unifying those that have appeared in the literature. This is achieved in Section 6.

In section 5 we study the case of flows from the trivial transformation semigroup on a 1-element set with the identity transformation. We will show that this is equivalent to a $\GM$ semigroup $S$ dividing the direct product of a group and its its right letter mapping image $\RLM(S)$. This case is intimately related to the lower bounds proved by Rhodes and Tilson in \cite{lowerbounds1, lowerbounds2}. The most important class of semigroups that have 1-point flows is the class of  inverse semigroups. Indeed, flows involve computing the minimal injective congruence on a transformation semigroup and inverse semigroups are faithfully represented by injective functions. We emphasize that many of the important results of inverse semigroup theory have a natural interpretation in flow theory. For a $\GM$ inverse semigroup $S$, the McAlister covering theorem that ensures that every inverse semigroup has an $E$-unitary cover \cite{McAlisterPstuff} is equivalent to the existence of a 1-point flow on $S$. Moreover $\RLM(S)$ is the fundamental image of $S$. The important decomposition theorem \cite{McAlReilly} is the restriction of the Presentation Lemma Decomposition Theorem \ref{PLFDT} to inverse semigroup theory. It is not surprising that important ideas and results in one field have been studied in neighboring fields. We wanted to give the explicit connections in this paper.

Section 7 and Section 8 of the paper are devoted to further examples of computing flows. In Section 7 we recast Tilson's 2-$\J$-class Theorem \cite{2J} in the language of flows. Tilson proved in the language of flows that a $\GM$  monoid $M$ with two non-zero $\J$-classes has complexity 1 if and only if it has a flow over the monoid $RZ(k)^{1}$ consisting of $k$ right zeros and an identity element, where $k$ is the number of orbits of the group of units of $M$ in its right action on the set of $\mathcal{L}$-classes of the 0-minimal $\J$-class of $M$. We then give some examples of calculating flows for this class of monoids. In Section 8, we look at examples of semigroups built from the character table of finite Abelian groups that also have flows over semigroups of the form $RZ(k)^{1}$. 

Computing flows are done in the Evaluation Transformation Semigroup defined in \cite{Trans} and uses many tools and definitions from this paper. In order to keep the ``flow" of the paper reasonable, we summarize the appropriate material from \cite{Trans} in the Appendix to the paper for the benefit of the reader. Many more examples of computing flows on semigroups or showing that they do not exist are contained in \cite{MasterList} that we encourage the reader to look at after reading this paper.
 
\section{Preliminaries}

All semigroups in this paper are finite. We assume a background in semigroup theory. This includes Rees Theorem and the Green relations. See \cite{howiebook, CP, CP2} for a classical approach. For results specifically related to finite semigroup theory see 
\cite{Arbib, Eilenberg, Lallement, qtheory, TilsonXI, TilsonXII}. We recall in more detail the Krohn-Rhodes Theorem and complexity of finite semigroups as these play a central role in this paper.

By a transformation semigroup (monoid) abbreviated ts (tm) we mean a pair $(Q,S)$ where $S$ is a subsemigroup (submonoid) of the monoid of all partial transformations $PT_{Q}$ acting on the right of $Q$. If $(P,T)$ and $(Q,S)$ are ts (tm), then their wreath product
$(P,T) \wr (Q,S)$ is the ts (tm) $(P \times Q, W)$ where $W = \{(f,s)\mid s \in S, f:Dom(s)\rightarrow T\}$ with action $(p,q)(f,s)=(p(sf),qf)$ whenever $q \in Dom(s), p \in Dom(sf)$ and undefined otherwise. $Dom(s)$ is the domain of the partial function $s$. It is useful to think of 
$W$ as consisting of $|Q| \times |Q|$ row-monomial matrices with entries in $T \cup \{0\}$. The matrix corresponding to $(f,s)$ is the row-monomial matrix whose entry in position $(q,qs)$ is $sf$ if $qs$ is defined and $q \in Dom(f)$ and $0$ otherwise. One easily sees that composition in $W$ corresponds to multiplication of the corresponding matrices, where we identify ``undefined'' with $0$.

For the formal definition of division of ts and tm see \cite{Eilenberg, qtheory}. In the case that all partial functions in $S$ and $T$ are total functions, this definition reduces to saying that $(P,T)$ divides $(Q,S)$ if $(P,T)$ is an image of a sub-ts (sub-tm) of $(Q,S)$. In particular, identifying a semigroup with its right regular representation we have that $T$ divides $S$ if and only if $T$ is an image of a subsemigroup (submonoid) of $S$. If $(P,T)$ divides $(Q,S)$ we write $(P,T) \prec (Q,S)$. A semigroup is aperiodic if all its subgroups are trivial.

\bt [Krohn-Rhodes 1962]

Every finite ts (tm) $(Q,S)$ divides a wreath product of groups and aperiodic semigroups. One can choose the groups to be simple groups that divide $S$ and the aperiodic semigroups to be sub-ts of $RZ(2)^1$ the transformation semigroup consisting of two states, the two constant maps and the identity function.

\et

We are lead to the definition of Krohn-Rhodes complexity.

\bd[\cite{KRannals}]

The Krohn-Rhodes complexity $Sc$ of a finite semigroup $S$ is the least number of groups in any decomposition of $S$ as a divisor of wreath
products of groups and aperiodic semigroups.

\ed

We say complexity of a finite semigroup in place of Krohn-Rhodes complexity. If the complexity of a semigroup $S$ is $n$, we write $Sc =n$. 
It is known \cite[Chapter 4]{qtheory} that the full transformation semigroup on $k$-elements, $T_k$, has complexity $k-1$. Therefore if we define $V_n$ to be the class of all semigroups of complexity at most $n$, then $V_{n}$ is properly contained in $V_{n+1}$ for all $n \geq 0$. As mentioned in the introduction, the papers \cite{complexity1} and \cite{complexityn} prove that the complexity of a finite semigroup is computable.

\subsection{Semilocal Theory, Group Mapping and Right Letter Mapping Semigroups}

Green's relations and Rees Theorem give the local structure of semigroups. Global semigroup theory is involved with decomposition theorems like the Krohn-Rhodes theorem and related tools to understand them. Semilocal theory is involved with how a semigroup acts on its Green classes. The main example is the Sch\"{u}tzenberger representation. We review the theory here. See Chapters 7 and 8 of \cite{Arbib} and 
Chapter 4 of \cite{qtheory} for more detailed information.

\bd

\begin{enumerate}
    
\item{A semigroup $S$ is right (left) transitive if $S$ has a faithful transitive right (left) action on some set $X$.}

\item{A semigroup $S$ is bi-transitive if $S$ has a faithful transitive representation on the right of some set $X$ and $S$ has a faithful transitive representation on the left of some set $Y$. We emphasize that we do not assume that $X=Y$.}

\end{enumerate}
\ed

\brm

In the literature these are called Right Mapping (denoted $\RM$), Left Mapping ($\LM$) and Generalized Group Mapping ($\GGM$) semigroups respectively. We will use these names from now on in this paper.

\erm

The most important example of a right transitive action is the Sch\"{u}tzenberger representation on an $\mathcal{R}$-class $R$ of $S$ whose definition we recall.

\medskip
For $r \in R, s \in S$, define 
    
		$$r\cdot s = \begin{cases}
      rs, & \text{if}\ rs \in R \\
      \text{undefined} & \text{otherwise}
    \end{cases}$$

There is the obvious dual notion of the left Sch\"{u}tzenberger representation of a semigroup on an $\mathcal{L}$-class.

\bl
 $S$ is an $\RM$, ($\LM$, $\GGM$) semigroup if and only if $S$ has a unique $0$-minimal regular ideal $I(S) \approx M^{0}(G,A,B,C)$ such that the right
 (left, right and left) Sch\"{u}tzenberger representation of $S$ is faithful on any $\mathcal{R}( \mathcal{L}, \mathcal{R} {\text{ and }} \mathcal{L})$-class $R (L, R {\text{ and }} L)$ class in $I(S) -\{0\}$.  
\el

\bd

A $\GGM$ semigroup is called {\em Group Mapping} written $\GM$ if the group $G$ in $I(S)$ is non-trivial.

\ed

Let $S$ be a $\GM$ semigroup. We identify a fixed $\mathcal{R}$-class in $I(S)$ with the set $G \times B$. We obtain a faithful transformation semigroup $(G \times B,S)$.

The action of $S$ on $G \times B$ induces an action of $S$ on $B$.  The faithful image of this action is called the 
{\em Right Letter Mapping} (written $\RLM(S)$) image of $S$.  We thus obtain a faithful transformation semigroup $(B,\RLM(S))$. We remark that $\RLM(S)$ has a 0-minimal regular ideal that is aperiodic and is the image of $I(S)$ under the morphism from $S$ to $\RLM(S)$. In particular, for any $\GM$ semigroup $S$, we have $|\RLM(S)| < |S|$. This allows for inductive proofs based on cardinality.

\bl

Let $S$ be a $GM$ semigroup. Then the following hold.

\begin{enumerate}

\item{$(G \times B,S)$ is a sub-ts of $G \wr (B,\RLM(S))$. Therefore, $Sc \leq 1+ \RLM(S)c$. }


\item{For every semigroup $T$, there is a $GM$ quotient semigroup $S$ such that $Tc=Sc$.} 

\end{enumerate}

\el

By induction on cardinality, we are lead to the following reduction theorem.

\bt
The question of decidability of complexity can be reduced to the case that $S$ is a $\GM$ semigroup and whether $Sc=\RLM(S)c$ or 
$Sc = 1+\RLM(S)c$.

\et

\brm

It is known that if $S$ is a completely regular semigroup, that is, a semigroup that is a union of its subgroups, then if $S$ is $\GM$ then it is always true that $Sc = 1+\RLM(S)c$ and there are completely regular semigroups of arbitrary complexity\cite[Chapter 9]{Arbib}. On the other hand, if $S$ is an inverse $\GM$ semigroup such that $\RLM(S)$ is not aperiodic, then $Sc=\RLM(S)c = 1$. The semigroup of all row and column monomial matrices over a non-trivial group $G$ is an example of such an inverse semigroup. We discuss the case of $\GM$ inverse semigroups in great detail later in this paper.

\erm

Flows are designed to give a necessary and sufficient, not a priori computable condition  for a $\GM$ semigroup $S$ to satisfy
$Sc=\RLM(S)c$.

\subsection{The Type II Subsemigroup and the Tilson Congruence}\label{redux.sec}

In this subsection we review the type II subsemigroup $S_{II}$ of a semigroup. It plays an important role in finite semigroup theory and a central role in flow theory. It was first defined in \cite{lowerbounds2}. In that paper it was proved that it is decidable if a regular element of a semigroup $S$ belongs to $S_{II}$. In particular, it followed that if $S$ is a regular semigroup, then membership in $S_{II}$ is decidable. Later 
Ash \cite{Ash} proved that membership in $S_{II}$ is decidable for all finite semigroups $S$.

The type II subsemigroup $S_{II}$ of $S$ is the smallest subsemigroup of $S$ containing all idempotents and closed under weak conjugation: if $xyx = x$ for $x,y \in S$, then $xS_{II}y \cup yS_{II}x \subseteq S_{II}$. Membership in $S_{II}$ is clearly decidable from this definition. Its importance stems from the following Theorem, where we give its original, not a priori decidable definition.

\bt

Let $S$ be a finite semigroup. Then $S_{II}$ is the intersection of all subsemigroups of the form $1\phi^{-1}$ where $\phi:S \rightarrow G$
is a relational morphism from $S$ to a group $G$.

\et

The aforementioned decidability results in \cite{lowerbounds2, Ash} prove that the two definitions we give define the same subsemigroup of a semigroup $S$. Here are some important properties of the type II subsemigroup. See \cite{qtheory} for proofs.

\bt

\begin{enumerate}

\item{Let $f:S \rightarrow T$ be a morphism between semigroups $S$ and $T$. If $s \in S_{II}$, then $sf \in T_{II}$. Therefore the restriction of $f$ to $S_{II}$ defines a morphism $f_{II}:S_{II} \rightarrow T_{II}$. This defines a functor on the category of finite semigroups.}

\item{If $S$ divides $T$, then $S_{II}$ divides $T_{II}$.}

\item{Assume that a semigroup $S$ divides a semidirect product $T*G$ where $T$ is a semigroup and $G$ is a group. Then $S_{II}$ divides $T$.}

\end{enumerate}

\et




A {\em congruence} on a transformation semigroup $(Q,S)$ is an equivalence relation $\approx$ on $Q$ such that if $q\approx q'$ and for $s \in S$, and both $qs$ and $q's$ are defined, then $qs\approx q's$. Every $s \in S$ defines a partial function on $\faktor{Q}{\approx}$ by 
$[q]_{\approx}s=[q's]_{\approx}$ if  $q's$ is defined for some $q' \in [q]_{\approx}$. The quotient $\faktor{(Q,S)}{\approx}$ has states $\faktor{Q}{\approx}$ and semigroup $T$ the semigroup generated by the action of all $s \in S$ on $\faktor{Q}{\approx}$. We remark that $T$ is not necessarily a quotient semigroup of $S$, but is in the case that $(Q,S)$ is a transformation semigroup of total functions. 

A congruence $\approx$ is called {\em injective} if every $s \in S$ defines a partial 1-1 function on $\faktor{Q}{\approx}$. It is easy to see that the intersection of injective congruences is injective. Therefore, there is a unique minimal injective congruence $\tau$ on any transformation semigroup $(Q,S)$. We call $\tau$ the Tilson congruence on a transformation semigroup because of the following proved in \cite{Redux}. It is central to the theory of flows.

\bt\label{redux}

Let $(G \times B,S)$ be a $\GM$ transformation semigroup. Then the minimal injective congruence $\tau$ on $S$ is defined as follows: 
$(g,b) \tau (g',b')$ if and only if there are elements $s,t \in S_{II}$ such that $(g,b)s=(g',b')$ and $(g',b')t=(g,b)$.

\et
 
\brm

\item{We can state this theorem by $(g,b)S_{II}=(g',b')S_{II}$. That is, $(g,b)$ and $(g',b')$ define the same ``right coset'' of $S_{II}$ on $G \times B$.}

\item{The proof in \cite{Redux} works on an arbitrary transitive transformation semigroup. We stated it in the case of $\GM$ transformation semigroups because that's how we will use it in this paper.}

\item{Theorem \ref{redux} can be used to greatly simplify the proof in \cite{lowerbounds2} for decidability of membership in $S_{II}$ for regular elements of an arbitrary semigroup.}

\erm

For later use we record the following corollary to Theorem \ref{redux}.

\bc

Let $(G \times B,S)$ be a $\GM$ transformation semigroup. Then $(g,b) \tau (g',b')$ if and only if there are elements 
$s,t \in S_{II} \cap I(S)$ such that $(g,b)s=(g',b')$ and $(g',b')t=(g,b)$. That is, we can choose the elements $s$ and $t$ in Theorem \ref{redux} to be in the 0-minimal ideal of $S$.

\ec

\proof The condition is sufficient by Theorem \ref{redux}. Conversely, assume that $(g,b)\tau (g',b')$. Then there are elements $s,t \in S_{II}$ such that $(g,b)s=(g',b')$ and $(g',b')t=(g,b)$. Since $I(S)$ is a 0-simple semigroup, there are idempotents $e,f \in I(S)$ such that 
$(g,b)e=(g,b)$ and $(g',b')f=(g',b')$. Therefore, $(g,b)es=(g,b)s=(g',b')$ and similarly $(g',b')ft =(g,b)$. Since $e,f$ are idempotents we have $es,ft \in S_{II}$. As $I(S)$ is a 0-minimal ideal, we also have $es,ft \in I(S)$. \Qed

%
%
%



\section{The Connection Between Rhodes Lattices and  Set-Partition Lattices} \label{subsection.SPRhodes}

Let $S$ be a $\GM$ semigroup with $0$-minimal ideal $M^{0}(G,A,B,C)$. The definition of flow gives a map from the state set of an automaton to a lattice associated with $S$. In the literature there are two such lattices. In this section we define these lattices and give the connections between them and show that they give equivalent notions of flow. The first is the set-partition lattice $SP(G \times B)$. This is the lattice whose elements are all pairs $(Y,\Pi)$ where $Y$ is a subset of $G\times B$ and $\Pi$ is a partition on $Y$. Here $(Y,\Pi) \leqslant (Z,\Theta)$ if $Y \subseteq Z$ and for all $y \in Y$, the $\Pi$ class of $y$ is contained in the $\Theta$ class of $y$. 

The second lattice is the Rhodes lattice $Rh_{B}(G)$. We review the basics. For more details see \cite{AmigoDowling}. Let $G$ be a finite group and $B$ a finite set. A partial partition on $B$ is a partition $\Pi$ on a subset $I$ of $B$. We also consider the collection of all functions
$F(B,G)$, $f:I \rightarrow G$ from subsets $I$ of $B$ to $G$. The group $G$ acts on the left of $F(B,G)$ by $(gf)(b) = gf(b)$ for $f \in F(B,G), g \in G, b \in Dom(f)$. An element $\{gf\mid f:I \rightarrow G, I \subseteq X, g \in G\}$ of the quotient set $F(B,G)/G$ is called a cross-section with domain $I$. It should be thought of as a projectification of a cross-section of the projection from $I$ to $B$ in the usual topological sense. An
SPC (Subset, Partition, Cross-section) over $G$ is a triple $(I,\Pi,f)$ where $I$ is a subset of $B$, $\Pi$ is a partition of $I$ and $f$ is a collection of cross-sections one for each $\Pi$-class $\pi$ with domain $\pi$. If the classes of $\Pi$ are $\{\pi_{1},\pi_{2}, \ldots, \pi_{k}\}$, then we sometimes write $\{(\pi_{1},f_{1}), \ldots, (\pi_{k},f_{k})\}$, where $f_{i}$ is the cross-section associated to $\pi_{i}$. For brevity we denote this set of cross-sections by $[f]_{\pi}$.  We let $Rh_{B}(G)$ denote the set of all SPCs on $B$ over the group $G$ union a new element $\Longrightarrow\Longleftarrow$  that we call {\em contradiction} and is the top element of the lattice structure on $Rh_{B}(G)$.  Contradiction occurs  because the join of two $SPCs$ need not exist. In this case we say that the contradiction is their join.

The partial order on $Rh_{B}(G)$ is defined as follows.

$(I,\pi,[f]_{\pi}) \leqslant (J,\tau,[h]_{\tau})$ if:
\begin{enumerate}
\item
{$I \subseteq J$} 
\item
{Every block of $\pi$ is contained in a (necessarily unique) block of $\tau$} 

\item{if the $\pi$-class $\pi_{i}$ is a subset of the $\tau$-class $\tau_{j}$, then the restriction of $h$ to $\pi_{i}$ equals $f$ restricted to $\pi_{i}$ as elements of $F(B,G)/G$. That is, $[h|_{\pi_{i}}] = [f|_{\pi_{i}}]$}
\end{enumerate}

See \cite[Section 3]{AmigoDowling} for the definition of the lattice structure on $Rh_{B}(G)$. The underlying set of the Rhodes lattice $Rh_{B}(G)$ minus the contradiction is the set underlying the Dowling lattice on the same set and group. The Dowling lattice has a different partial order. For the connection between Rhodes lattices and Dowling lattices see \cite{AmigoDowling}.

We note that $SP(G \times B)$ is isomorphic to the Rhodes lattice $Rh_{G\times B}(1)$ of the trivial group over the set $G \times B$. We need only note that a cross-section to the trivial group is a partial constant function to the identity and can be omitted, leaving us with a set-partition pair. There are no contradictions for Rhodes lattices over the trivial group, and in this case the top element is the pair $(G \times B,(G \times B)^{2})$. Despite this, we prefer to use the notation $SP(G \times B)$ instead of $Rh_{G\times B}(1)$.

Conversely, we can find a copy of the meet-semilattice of $Rh_{B}(G)$ as a meet subsemilattice of $SP(G\times B)$. We begin with the following important definition.

\bd
A subset $X$ of $G \times B$ is a cross-section if whenever $(g,b),(h,b) \in X$ then $g = h$. That is, $X$ defines a cross-section of the projection $\theta:X \rightarrow B$. Equivalently, $X^{\rho} \subseteq B \times G$, the reverse of $X$, is the graph of a partial function $f_{X}:B \rightarrow G$.  An  element $(Y,\Pi) \in SP(G \times B)$ is a cross-section if  every partition class $\pi$ of $\Pi$ is a cross-section.  
\ed

From the semigroup point of view, a cross-section is a partial transversal of the $\mathcal{H}$-classes of the distinguished $\mathcal{R}$-class, $R = G \times B$ of a $\GM$ semigroup. That is, the $\mathcal{H}$-classes of $R$ are indexed by $B$ and a cross-section picks at most one element from each $\mathcal{H}$-class. 

An element $(Y,\Pi)$ that is not a cross-section is called a {\em contradiction}. That is, $(Y,\Pi)$ is a contradiction if some $\Pi$-class $\pi$ contains two elements $(g,b),(h,b)$ with $g \neq h$. We note that the set of cross-sections is a meet subsemilattice of $SP(G \times B)$ and the set of contradictions is a join subsemilattice of $SP(G \times B)$.

We will identify $B$ with the subset $\{(1,b)\mid b \in B\}$ of $G \times B$, which as above we think of as a system of representatives of the $\mathcal{H}$-classes of the distinguished $\mathcal{R}$-class. Note then that $G\times B$ is the free left $G$-act on 
$B$ under the action $g(h,b)=(gh,b)$. This action extends to subsets and partitions of $G\times B$ . Thus $SP(G\times B)$ is a left $G$-act. An element $(Y,\Pi)$ is {\em invariant} if $G(Y,\Pi) = (Y,\Pi)$. It is easy to see that $(Y,\Pi)$ is invariant if and only if:

\begin{enumerate}
 \item{$Y=G \times B'$ for some subset $B'$ of $B$.}  

 \item{For each $\Pi$-class $\pi$, $G\pi \subseteq \Pi$ and is a partition of $G\times B''$ where $B''\subseteq B'$.}
\end{enumerate}

Thus $(Y,\Pi)$ is invariant if and only if $Y=G\times B'$ for some subset $B'$ of $B$ and there is a partition $B_{1},\ldots B_{n}$ of $B'$ such that for all $\Pi$ classes $\pi$, $G\pi$ is a partition of $G\times B_{i}$ for a unique $1 \leqslant i \leq n$. 

Let $CS(G\times B)$ be the set of invariant cross-sections $(Y,\Pi)$ in $SP(G \times B)$. $CS(G\times B)$ is a meet subsemilattice of $SP(G\times B)$. We claim it is isomorphic to the meet-semilattice of $Rh_{B}(G)$. Indeed, let $(G\times B',\Pi)$ be an invariant cross-section. Using the notation above, pick  $\Pi$-classes, $\pi_{1}, \ldots , \pi_{n}$ such that $G\pi_{i}$ is a partition of $G\times B_{i}$ for $i=1, \ldots n$. Since $\pi_{i}$ is a cross-section its reverse is the graph of a function $f_{i}:B_{i}\rightarrow G$. We map 
$(G\times B',\Pi)$ to $(B',\{B_{1},\ldots , B_{n}\},[f]_{\{B_{1},\ldots , B_{n}\}}) \in Rh_{B}(G)$ where the component of $f$ on $B_{i}$ is $f_i$. It is clear that this does not depend on the representatives $\pi$ and we have a well-defined function from $F:CS(G \times B) \rightarrow Rh_{B}(G)$. It is easy to see that $F$ is a morphism of meet-semilattices.

Conversely, let $(B',\Theta,[f]) \in Rh_{B}(G)$. We map $(B',\Theta,[f])$ to $(G \times B',\Pi) \in CS(G \times B)$ where $\Pi$ is defined as follows. If $\theta$ is a $\Theta$-class, let $\widehat{\theta} =\{(bf,b)\mid  b \in \theta\}$. Let $\Pi$ be the collection of all subsets of the form $g\widehat{\theta}$ where $\theta$ is a $\Theta$-class and $g \in G$. The proof that $(G \times B',\Pi) \in CS(G \times B)$ and that this assignment is the inverse to $F$ above and gives an isomorphism between the meet-semilattice of $Rh_{B}(G)$ and $CS(G \times B)$ is straightforward and left to the reader. Furthermore, if the join of two $SPC$ in $Rh_{B}(G)$ is an $SPC$ (that is it is not the contradiction) then this assignment preserves joins. Therefore we can identify the lattice $Rh_{B}(G)$ with $CS(G \times B)$ with a new element $\Longrightarrow\Longleftarrow$ added when we define the join of two elements of $CS(G \times B)$ to be $\Longrightarrow\Longleftarrow$ if their join in $SP(G\times B)$ is a contradiction. We use $CS(G\times B)$ for this lattice as well.

We record this discussion in the following Proposition.

\bp\label{Updown}
The set-partition lattice $SP(G\times B)$ is isomorphic to the Rhodes lattice $Rh_{G\times B}(1)$. The Rhodes lattice $Rh_{B}(G)$ is isomorphic to the lattice $CS(G\times B)$.
\ep

\section{The Presentation Lemma, Flows and the Flow Decomposition Theorem}

We review the definition of flows from  an automaton with alphabet $X$  to the set-partition $SP(G\times B)$ and to the Rhodes lattice $Rh_{B}(G)$. For more details, see \cite[Sections 2-3]{Trans}.

Let $(G \times B,S)$ be the transformation semigroup associated to a $\GM$ semigroup $S$ and let $X$ be a generating set for $S$. By a deterministic automaton we mean an automaton such that each letter defines a partial function on the state set.

\bd

Let $\mathcal{A}$ be a deterministic automaton with state set $Q$ and alphabet $X$. 
A {\em flow} to the lattice $SP(G\times B)$ on $\mathcal{A}$ is a function $f:Q \rightarrow SP(G \times B)$ such that for each $q \in Q, x \in X$, with $qf =(Y,\pi)$ and $(qx)f = (Z,\theta)$ we have:

\begin{enumerate}

\item{For all $(g,b) \in G \times B$, there is a $q \in Q$ such that $(g,b) \in qf$.}
 
\item{$Yx \subseteq Z$.}

\item{Multiplication by $x$ considered as an element of $S$ induces a partial 1-1 map from $Y/\pi$ to $Z/\theta$.

That is, for all $y,y' \in Y$ we have $y,y'$ are in a $\pi$-class if and only if $yx,y'x$ are in a $\theta$-class whenever $yx,y'x$ are both defined.}

\item{{\bf The Cross Section Condition:} For all $q \in Q$, $qf$ is a cross-section.}

\end{enumerate}
\ed

We use Proposition \ref{Updown} to give a definition of a flow to the Rhodes lattice $Rh_{B}(G)$.

\bd

Let $\mathcal{A}$ be a deterministic automaton with state set $Q$ and alphabet $X$. 
A {\em flow} to the lattice $Rh_{B}(G)$ is a flow to $SP(G\times B)$ such that for each state $q$, $qf \in CS(G\times B)$. That is, $qf$ is a $G$-invariant cross-section.
   
\ed

\brm

    It follows from the original statement of the Presentation Lemma \cite{AHNR.1995} that if $S$ has a flow with respect to some automaton over $SP(G\times B)$, then it has a flow from the same automaton over $Rh_{B}(G)$. The proof of the equivalence of the Presentation Lemma and Flows in Section 3 of \cite{Trans} works as well for the Presentation Lemma in the sense of \cite{AHNR.1995}. 

\erm

We will use all the terminology and concepts from \cite{Trans}. If $\mathcal{A}$ is an automaton with alphabet $X$ and state set $Q$, recall that its completion is the automaton $\mathcal{A}^{\square}$ that adds a sink state $\square$ to $Q$ and declares that $qx = \square$ if $qx$ is not defined in $\mathcal{A}$. A flow on $\mathcal{A}$ is a complete flow if on $\mathcal{A}^{\square}$ extending $f$ by letting $\square f = (\emptyset, \emptyset)$ the bottom of both the set-partition and Rhodes lattices remains a flow and furthermore, for all $(g,b) \in G \times B$, there is a $q\in Q$ such that $(g,b) \in Y$, where $qf=(Y,\Pi)$. All flows in this paper will be complete flows.

 Flows are related to the Presentation Lemma \cite{AHNR.1995}, \cite[Section 4.14]{qtheory}. They give a necessary and sufficient condition for $Sc =\RLM(S)c$ where $S$ is a $\GM$ semigroup.

\bt\label{PLflow}[The Presentation Lemma-Flow Version]

Let $(G \times B,S)$ be a $GM$ transformation semigroup with $S$ generated by $X$. Let $k> 0$ and assume that $\RLM(S)c = k$.  Then $Sc = k$ if and only if there is an $X$ automaton $\mathcal{A}$ whose transformation semigroup $T$ has complexity strictly less than $k$ and a complete flow $f:Q \rightarrow Rh_{B}(G)$.


\et


The proof of Theorem \ref{PLflow} uses the decomposition theorem (Theorem \ref{PLFDT} below), that is of independent interest and central to the results in this paper. The proof was originally given in \cite{AHNR.1995}. The decomposition can be more easily proved using Steinberg's Slice Theorem \cite{Slice2}, \cite[Section 2.7]{qtheory}. See also Exercises 4.14.22-4.14.24 of \cite{qtheory} which give the connection between the Presentation Lemma and the Slice Theorem. We make these connections precise in Section \ref{unified}. 

The Slice Theorem gives a necessary and sufficient condition for a semigroup to divide a direct product. Recall that a division $d:S
 \rightarrow T$ is an injective relational morphism. This means that if $sd \cap s'd \neq \emptyset$ then $s=s'$. Equivalently, $d$ is the inverse of a surjective morphism from a subsemigroup of $T$
 onto $S$, so that $S$ divides $T$ in the usual sense of semigroup theory if and only if there is a division $d:S \rightarrow T$. See 
\cite{qtheory} for details. Recall that if $f:S \rightarrow T$ is a relational morphism, then $graph(f)=\{(s,t) \mid t \in sf\}$ is a subsemigroup of $S \times T$.


\bt{\bf The Slice Theorem}\label{benslice}

Let $f:S \rightarrow T$ be a relational morphism. Then $f=d\pi$ where $d:S \rightarrow U \times T$ is a division and $\pi:U \times T \rightarrow T$ the projection to some semigroup $U$ if and only if there is a relational morphism $g:graph(f) \rightarrow U$ satisfying the Slice Condition:
$$(s,t)g \cap (s',t)g \neq \emptyset {\text{ implies }} s=s'.$$

\et

See \cite[Section2.7]{qtheory} for a proof of the Slice Theorem. If $f$ in Theorem \ref{benslice} is a functional morphism, then we get the following Corollary. We just note that if $f:S \rightarrow T$ is a functional morphism then $graph(f)$ is isomorphic to $S$ via the projection $graph(f) \rightarrow S$. That is, we identify $(s,sf)$ with $s$.

\bc \label{fslice}

Let $f:S \rightarrow T$ be a functional morphism. Then $f=d\pi$ with $d:S \rightarrow U \times T$ is a division and $\pi:U \times T \rightarrow T$ the projection for some semigroup $U$ if and only if there is a relational morphism $g:S \rightarrow U$ satisfying the Slice Condition: $sg \cap s'g \neq \emptyset$ and $sf=s'f$ implies $s=s'$.
\ec

We reformulate in the case the morphism in Corollary \ref{fslice} is surjective.

\bc\label{slice}

Let $S$ and $U$ be semigroups,  $\equiv$ a congruence on $S$ and $\phi:S \rightarrow U$ a relational morphism to $U$. Then 
$S \prec U \times \faktor{S}{\equiv}$ if and only if
for all $s,s' \in S$, if $s\equiv s'$ and $s\phi \cap s'\phi \neq \emptyset$ then $s=s'$. 

\ec

\proof Take $T = \faktor{S}{\equiv}$ in Corollary \ref{fslice} and $f:S \rightarrow \faktor{S}{\equiv}$ the natural morphism.


We apply the Slice Theorem to obtain the following Decomposition Theorem.

\bt

\label{PLdecomp}
 Let $(G \times B,S)$ be a $\GM$-transformation semigroup and let $k > 0$.  Then $Sc \leqslant k$ if and only if $\RLM(S)c \leqslant k$ and there is a transformation semigroup $T$ with $Tc \leqslant k-1$ and a group $H$ such that 
$$ (G \times B,S) \prec (H \wr T) \times (B,\RLM(S)) $$
\et

\proof If $\RLM(S)c \leqslant k$ and there is a transformation semigroup $T$ with $Tc \leqslant k-1$ such that 
$ (G \times B,S) \prec (G \wr Sym(B) \wr T) \times (B,\RLM(S))$, then $Sc \leqslant \max\{(G \wr Sym(B) \wr T)c,\RLM(S)c\} \leq k$ by the assumptions on $\operatorname{\RLM}$ and $T$.

Conversely, let $(G \times B,S)$ be a $\GM$-transformation semigroup with $S$ generated by $X$ and $Sc \leqslant k$. Then 
$S \prec A \wr H \wr T$ where $A$ is an aperiodic ts, $H$ is a transformation group and $Tc \leqslant k-1$. By considering the projection $\pi:A \wr H \wr T \rightarrow H \wr T$ we obtain an aperiodic relational morphism of ts $\phi:(G \times B,S) \rightarrow H \wr T$. We also consider the congruence $\equiv$ of the morphism from $S$ to $\RLM(S)$.

We verify the conditions in Corollary \ref{slice}. Assume that $s \equiv s'$ and that $s\phi \cap s'\phi \neq \emptyset$. Let $(g,b) \in G \times B$. Since $s \equiv s'$, we have $(g,b)s$ is defined if and only if $(g,b)s'$ is defined. Assume both are defined. Then 
$(g,b)s=(g(bs),bs)$ and $(g,b)s'=(g(bs'),bs')$. But $bs =bs'$ since $s\equiv s'$. Therefore, $(g,b)s \mathcal{H} (g,b)s'$. Let 
$t \in s\phi \cap s'\phi$ and $u \in (g,b)\phi$. Then $ut \in (g,b)s\phi \cap (g,b)s'\phi \neq \emptyset$. Since $\phi$ is an aperiodic relational morphism, it follows that $(g,b)s=(g,b)s'$, since $\phi$ separates elements in regular $\mathcal{H}$-classes. Therefore, $s$ and $s'$
act the same on $G \times B$ and thus $s=s'$ since $S$ is $\GM$. We have shown that the conditions of Lemma \ref{slice} hold and thus
$ (G \times B,S) \prec (H \wr T) \times (B,\RLM(S))$. \Qed

Flows allow us to show that we can take the group $H$ in Theorem \ref{PLdecomp} to be $G \wr Sym(B)$. In particular the size of the group in Theorem \ref{PLdecomp} is bounded as a function of the input $(G \times B,S)$.

\bt \label{GwrSymB}

Let $(G \times B,S)$ be a $\GM$ transformation semigroup generated by $X$. Let $(Q,T)$ be a transformation semigroup with $T$ also generated by $X$. Then there is a flow $\Phi:Q \rightarrow Rh_{B}(G)$ if and only if $(G \times B, S) \prec (G \wr Sym(B) \wr T) \times \RLM(S)$.
 
\et

\proof Let $\Phi:Q \rightarrow Rh_{B}(G)$ be a flow. For each $q \in Q$, let $qF= (B_{q},\Pi_{q},F_{q})$ where $B_{q} \subseteq B$,
$\Pi_{q}$ is a partition of $B_{q}$ with classes $\pi_{q}^{1}, \ldots , \pi_{q}^{n_{q}}$ and 
$F_{q}=\{f_{q}^{i}:\pi_{q}^{i}\rightarrow G \mid i = 1 \ldots , n_{q}\}$ are representatives for the cross-section of $qF$. Pick a representative 
$b_{q}^{i} \in \pi_{q}^{i}, i = 1, \ldots , n_{q}$. 

Let $x \in X$. Then by the definition of a flow, there is a partial 1-1 function 
$\faktor{B_{q}}{\Pi_{q}}\rightarrow \faktor{B_{qx}}{\Pi_{qx}}$ that preserves cross-sections. This means that if $\pi_{q}^{i}x \neq \emptyset$, then there are a unique $j$ and a unique $g_{i,j} \in G$ such that $\pi_{q}^{i}x \subseteq \pi_{qx}^{j}$ and $f_{q}^{i}x=g_{i,j}f_{qx}^{j}$.
This information can be encoded in an $n_{q} \times n_{qx}$ row and column monomial matrix $M_{q,qx}$ whose $(i,j)$ entry is $g_{i,j}$. Since both 
$n_{q}$ and $n_{qx}$ are at most the cardinality of $B$, we can embed $M_{q,qx}$ into a monomial matrix 
$\mathcal{M}_{q,qx} \in G \wr Sym(B)$. Let $x \in X$. Then there are elements $s_{x} \in S$ and  $t_{x} \in T$ such that $t_{x}$ covers $s_{x}$. Let $\Theta:S \rightarrow G \wr Sym(B) \wr T$ be the relational morphism generated by sending $s_{x}$ to $(\rho_{x},t_{x})$ where $q\rho_{x} = \mathcal{M}_{q,qx}$. We claim that $\Theta$ satisfies the conditions of Corollary 
\ref{slice} relative to the congruence $\equiv$ induced by the morphism from $S$ to $RLM(S)$. This claim and Corollary \ref{slice} imply that
$(G \times B, S) \prec (G \wr Sym(B) \wr T) \times \RLM(S)$.

Assume then that $s \equiv s'$ and that $s\Theta \cap s'\Theta \neq \emptyset$. Suppose that $(\rho,t)$ covers both $s$ and $s'$. Let 
$(g,b,q)$ be a state in the graph of $\Theta$. Assume that $qt$ is defined. Then $(g,b)s =(g(bs),bs), (g,b)s'=(g(bs'),bs')$. Since $s \equiv s'$, we have $bs =bs'$. Therefore, $(g,b)s\mathcal{H}(g,b)s'$. Since $\Phi$ is a flow, it follows that $(g,b)s=(g,b)s'$.  Since every $(g,b)$ in $q\Phi$ for some $q \in Q$, we have $(g,b)s=(g,b)s'$ for all $(g,b) \in G \times B$ and it follows that $s=s'$ since $S$ is $\GM$. \Qed






The next Theorem follows immediately from Theorem \ref{PLflow} and Theorem \ref{GwrSymB}.

\bt[{The Presentation Lemma Flow Decomposition Theorem}]
\label{PLFDT}
 Let $(G \times B,S)$ be a $\GM$-transformation semigroup and let $k > 0$.  Then $Sc \leqslant k$ if and only if $\RLM(S)c \leqslant k$ and there is a transformation semigroup $T$ with $Tc \leqslant k-1$ such that 
$$ (G \times B,S) \prec (G \wr Sym(B) \wr T) \times (B,\RLM(S)).$$
\et

\brm

We note that given $(G \times B,S)$, we can effectively compute $(B,\RLM(S))$ and $G \wr Sym(B)$. Therefore the decomposition $(G \times B,S) \prec (G \wr Sym(B) \wr T) \times (B,\RLM(S))$
from Theorem \ref{PLFDT} reduces decidability of complexity to finding an appropriate $T$. Flows give us the tools to find such a $T$.

\erm

From the decomposition $(G \times B,S) \prec (G \wr Sym(B) \wr T) \times (B,\RLM(S))$ we obtain a subsemigroup $V$ of $(G \wr Sym(B) \wr T) \times (B,\RLM(S))$ and a surjective morphism $\phi:V \rightarrow S$. Let $\pi$ be the projection from 
$(G \wr Sym(B) \wr T) \times (B,\RLM(S))$ to $(G \wr Sym(B) \wr T)$. The subsemigroup $Res(f) = V\pi$ of 
$(G \wr Sym(B) \wr T)$ is called the resolution semigroup of the flow. The minimal ideal $K(Res(f))$ is called the resolution ideal of $f$ and resolves the $0$ of $S$. We then have a relational morphism $(\phi)^{-1}\pi:S \rightarrow Res(f)$ called the resolving morphism. We'll look at examples of resolutions in the next section.

In particular, in the case of inverse semigroups, we will see that $S$ has a flow over the trivial semigroup and the semigroup $V$ is an $E$-unitary cover of $S$. $Res(f)=K(Res(f))$ is a group such that $S$ has a relational morphism $\Theta:S \rightarrow Res(f)$ with $1(\Theta)^{-1} = E(S)$, where $E(S)$ is the semilattice of idempotents of $S$. We detail this in the next section and more generally discuss the $\GM$ semigroups that have flows over the trivial semigroup.

\section{One-point flows}
In this section we analyze $\GM$ semigroups that have one-point flows. This means that there is a flow from the trivial automaton with one state and the identity transformation. These are precisely the semigroups that have flows $f$ whose resolution semigroup is a group. This is discussed in a different language in \cite{lowerbounds2}, \cite{Redux} and \cite{Mcalisterschutz} without the use of flows, which these papers pre-dated. 

An important class of semigroups that have 1-point flows is the class of inverse semigroups. This is essentially the theory of $E$-unitary covers of inverse semigroups. We discuss this from the flow theoretic perspective in detail in this section giving the connections to the classical theory of inverse semigroups.

Let $S$ be a $\GM$ semigroup with distinguished ideal $I(S)\approx M^{0}(G,A,B,C)$. We have the corresponding transformation semigroup $X=(G \times B,S)$. In Section \ref{redux.sec} we defined the Tilson congruence $\tau$ \cite{Redux} on $X$ defined by $(g,b) \tau (g',b')$ if and only if there are $s,t \in S_{II}$ such such that $(g,b)s = (g',b')$ and $(g',b')t=(g,b)$ and noted that we can choose 
$s,t \in S_{II} \cap I(S)$. We also note that if $(g,b)\tau (g', b')$, then $(hg,b)\tau (hg',b')$ for all $h \in G$. This follows immediately from the definition. It follows that $\faktor{X}{\tau}$ is a left $G$-set. We recall that the key property of $\tau$ is that it is the smallest injective congruence on $X$. The connection with flows is given by the following definition.


%

%
%
%
%

\bd

We call $\tau$ a cross-section if for all $g,h \in G, b \in B$, if $(g,b)\tau(h,b)$ then $g=h$. By the left invariance mentioned above, this is equivalent to $(1,b)\tau (g,b)$ implies $g=1$ for all $g \in G, b\in B$.

\ed

Thus $\tau$ is a cross-section if and only if each $\tau$-class $[(g,b)]_{\tau}$ is a cross-section in the topological sense of the projection from $[(g,b)]_{\tau}$ to $B$.

\bt \label{onepoint}

Let $S$ be a $\GM$ semigroup. Then the following conditions are equivalent.

\begin{enumerate}

\item{$S_{II} \cap I(S)$ is aperiodic.}

\item{$\tau$ is a cross section.}

\item{There is a one-point flow $\Phi:(\{1\},\{1\}) \rightarrow Rh_{B}(G)$.}

\end{enumerate}

\et

\proof ($1 \rightarrow 2$) Assume that $S_{II} \cap I(S)$ is aperiodic. If $(1,b) \tau (g,b)$, then by the remark above there are elements in 
$s,t \in S_{II} \cap I(S)$ such that $(1,b)s=(g,b)$ and $(g,b)t=(1,b)$. By Rees Theorem, we have $s$ is an element of a group in the 
$\mathcal{L}$ class of $(1,b)$. That is, $s=(a,h,b)$ for some $a \in A, h \in G$ with $C(b,a) \neq 0$. Since $S_{II} \cap I(S)$ is 
aperiodic, it follows that $s$ is an idempotent and thus a right identity for $(1,b)$. Therefore $g=1$ and $\tau$ is a cross section. 

($2 \rightarrow 3$) Assume that $\tau$ is a cross section. We first claim that $\tau$ induces an equivalence relation on $B$, denoted by 
$\tau_{B}$, by defining $b\tau_{B} b'$ if there are $g,h \in G$ such that $(g,b) \tau (h,b')$. It is clear that $\tau_{B}$ is reflexive and symmetric.
Assume that $b\tau_{B} b'$ and $b'\tau_{B}b''$. Then there are $g,h,k,l \in G$ such that $(g,b) \tau (h,b')$ and $(k,b') \tau (l,b'')$. By left $G$-invariance of $\tau$, we have $(h,b') \tau (hk^{-1}l,b'')$. Therefore, $(g,b) \tau (hk^{-1}l,b'')$ and thus $\tau_{B}$ is transitive.    

Define $\Phi:\{1\} \rightarrow Rh_{B}(G)$ by $1\Phi$ is the $SPC$ whose set is $B$ and whose partition is the partition defined by $\tau_{B}$. Let $\pi = \{b_{1},\ldots, b_{k}\}$ be a $\tau_{B}$ class. Define $f_{\pi}:\pi \rightarrow G$ by $b_{i}f_{\pi} = g_{i}$ if 
$(1,b_{1}) \tau (g_{i},b_{i})$. By left $G$-invariance of $\tau$ it follows that $f_{\pi}$ defines a well-defined cross section. This data defines the function $\Phi:\{1\} \rightarrow Rh_{B}(G)$. $\Phi$ is a flow because $\tau$ is the smallest injective congruence on $X$.

($3 \rightarrow 1$) Assume that $\Phi:(\{1\},\{1\}) \rightarrow Rh_{B}(G)$ is a flow. Thus, we have an SPC $1\Phi$ whose set is $B$. Define a relation
$\rho_{\Phi}$ on $G \times B$ by $(g,b)\rho_{\Phi} (h,b')$ if $\{b,b'\}$ belong to a partition class of $\pi$ of $1\Phi$ and the cross section $f_{\pi}$ with domain $\pi$ has value $g$ on $b$ and $h$ on $b'$. The proof that $\rho_{\Phi}$ is an 
equivalence relation is similar to the proof in the last paragraph using the left $G$-invariance of cross-sections. It follows easily from the assumption that $\Phi$ defines a flow that $\rho_{\Phi}$ is an injective congruence on $G \times B$. Therefore $\tau \subseteq \rho_{\Phi}$ by Theorem \ref{redux}. Clearly, if $(g,b)\rho_{\Phi}(h,b)$, then $g=h$. Therefore the same property holds for $\tau$. In particular, if $(1,b) \tau (g,b)$, then 
$g=1$. It follows that if $x = (a,h,b)$ is a group element in $S_{II} \cap I(S)$, then $x$ is an idempotent, since groups in $I(S)$ act freely on the right of their $\mathcal{L}$-class. Therefore, $S_{II} \cap I(S)$ is aperiodic and the proof is complete. \qed

\brm

We note another condition equivalent to the three in Theorem \ref{onepoint}: If $I(S) \cap S_{II}$ is aperiodic, then in particular 
$I(S) \cap IG(S)$ is aperiodic, where $IG(S)$ is the idempotent generated subsemigroup of $S$. This is well known to be equivalent to $I(S)$
 having a Rees matrix representation where all elements of its structure matrix are 0 or 1. Assuming this normalization and taking the 
 congruence $\tau$ on $G \times B$ every element of $S_{II}$ acts on $G \times B$ by a row monomial matrix with entries in $\{0,1\}$. This
 is equivalent to there being a flow in which all the cross-sections are the constant cross section having all entries equal to 1. 

\erm

The following Corollary follows from Theorem \ref{onepoint} and Theorem \ref{PLFDT}.

\bc

Let $S$ be a $\GM$ semigroup such that $S_{II} \cap I(S)$ is aperiodic. Then $S$ divides 
$G \wr Sym(B) \times \RLM(S)$. In particular, the resolution ideal of a 1-point flow is a subgroup of $(G \wr Sym(B))$.

\ec

We define another kind of triviality for a flow.

\bd

Let $X= (G \times B,S)$ be $\GM$ transformation semigroup. $X$ is called {\em $0$-constant} if every $s \in S$ acts on $G \times B$ by a row-monomial matrix such that each entry belongs to $\{0,g\}$ for some $g \in G$.
\ed

\bp

$(G \times B,S)$ is 0-constant if and only if $S \prec G \times \RLM(S)$.

\ep

\proof We identify $\RLM(S)$ with the monoid of $B \times B$ row-monomial matrices over $\{0,1\}$ that arise by replacing the non-zero entries of each row monomial matrix in $S$ by 1. The minimal ideal $I$ of $G \times \RLM(S)$ is $G \times \{0\}$. The map from 
$G \times \RLM(S)/I$ to the semigroup of all 0-constant matrices over $G$ by sending $(g,M)$ to $gM$, where $M$ is a matrix in $\RLM(S)$ is an isomorphism. The result follows easily from this identification. \Qed

\bcomp

We now give a general method for computing the resolution semigroup of a $\GM$ semigroup $S$ with a 1-point flow. Let $S$ be a $\GM$ semigroup with a 1-point flow. Let $R:S \rightarrow H$ be a surjective relational morphism to a group $H$ such that 
$1R^{-1} \cap I(S) = S_{II} \cap I(S)$ an aperiodic subsemigroup of $I(S)$ by Theorem \ref{onepoint}. We claim that $S$ divides $H \times \RLM(S)$. Let $f:S \rightarrow \RLM(S)$ be the canonical morphism. 

By the Slice Theorem \ref{benslice} we need to show that if $sR \cap s'R \neq \emptyset$ and $sf = s'f$ then $s=s'$.
Let $(g,b) \in G \times B$. Since $sf = s'f$, $(g,b)s$ is defined if and only if $(g,b)s'$ is defined.  Furthermore if  $(g,b)s=(g(b)s,bs), (g,b)s'=(g(b)s,bs')$ it follows that $bs=bs'$. Therefore $(g,b)s \mathcal{H} (g,b)s'$. Let $h \in sR \cap s'R$ and $k \in (g,b)R$. Then 
$hk \in  ((g,b)s))R \cap ((g,b)s'))R$. Therefore $(hk)R^{-1}$ contains both $(g,b)s$ and $(g,b)s'$. But since $1R^{-1} \cap I(S)$ is aperiodic it follows easily that for all $x \in H$, $xR^{-1}$ intersects any $\mathcal{H}$ class of $I(S)$ in at most one point. Therefore $(g,b)s=(g,b)s'$. Since $(g,b)$ is an arbitrary element of $G \times B$ and $S$ is $\GM$ it follows that $s=s'$. This completes the proof that $S$ divides $H \times \RLM(S)$. Therefore $H$ is the resolution ideal of this decomposition. Let $T =\{(s,g)\mid (s,g) \in R\}$ be the graph of $R$. Then the projections to both $S$ and $H$ are surjective. $T$ is then the resolution semigroup of $S$.

\ecomp

\bcomp

We now show how a 1-point flow gives a resolution. We call this the resolution associated to the flow. Let then $X = (G \times B,S)$ be a $\GM$ transformation semigroup. Assume that $X$ has a 1-point flow. Start by computing the congruence $\tau$. By Theorem \ref{onepoint}, the classes of $\tau$ are a cross-section. Since $\tau$ is $G$-invariant, it induces an equivalence relation on $B$ that we also denote by $\tau$. Furthermore, $\faktor{X}{\tau}$ is an injective transformation semigroup with state set $\faktor{G \times B}{\tau}$. Thus each $s \in S$ defines a partial 1-1 map $\bar{s}$ on $\faktor{G \times B}{\tau}$. Since $s$ can be represented by a $B \times B$ row monomial matrix over $G$, $\bar{s}$ can be represented by a $\faktor{B}{\tau} \times \faktor{B}{\tau}$ row {\em and} column monomial matrix over $G$. 

We let $\rho(s)$ to be the set of all $\faktor{B}{\tau} \times \faktor{B}{\tau}$ monomial matrices (all rows and columns have exactly one non-zero entry) over $G$ that contain $\bar{s}$ as a submatrix. It is not difficult to show that $\rho:S \rightarrow G \wr Sym(\faktor{B}{\tau})$ is a relational morphism where we identify $G \wr Sym(\faktor{B}{\tau})$ with the group of $\faktor{B}{\tau} \times \faktor{B}{\tau}$ monomial matrices  over $G$. We claim that $S\rho$ is a resolution group. Indeed, $\tau$ being a cross-section immediately gives the conditions of the Slice Theorem relative to the morphism from $S$ to $\RLM(S)$. We suggest that the reader use this method to compute resolution semigroups for the examples of semigroups with 1-point flows that we give below.

\ecomp

The following follows from the computation above of the resolution associated to a flow. It is a refinement of Theorem \ref{PLFDT}.

\bt

Let $(G \times B,S)$ be $\GM$ with a 1-point flow. Then 
$$S \prec (G \wr Sym(\faktor{B}{\tau})) \times \RLM(S)$$

\et

The following Theorem gives equivalences for being in the pseudovariety $Ap*Gp$ of all semigroups $S$ that divide a wreath product $A \wr G$ where $A$ is aperiodic and $G$ is a group. Proofs can be found in  \cite{lowerbounds2, complex1/2, qtheory}.

\bp\label{Ap*Gp}

Let $S$ be a finite semigroup. The following conditions are equivalent.

\begin{enumerate}

\item{$S \in Ap*Gp$.}

\item{There is a relational morphism $R:S \rightarrow G$ from $S$ to a group $G$ such that $1R^{-1}$ is aperiodic.}

\item{$S_{II}$  is aperiodic.}

\end{enumerate}

\ep

%

The following Corollary follows immediately from Theorem \ref{onepoint} and Proposition \ref{Ap*Gp}.
 
\bc

Let $S$ be a $\GM$ semigroup in $Ap*Gp$. Then $S$ divides $(G \wr Sym(B)) \times \RLM(S)$.

\ec

We now show that $Ap*Gp$ is the unique largest pseudovariety such that each of its members has a 1-point cover. 

\bp

Let $S$ be a semigroup that does not belong to $Ap*Gp$. Then some $\GM$ image of $S$ does not have a 1-point flow. Consequently, $Ap*Gp$ is the unique largest pseudovariety such that each of its $GM$ images has a 1-point flow.

\ep

\proof Since $S$ is not in $Ap*Gp$ there is some $\mathcal{J}$-class $J$ of $S$ such that $J \cap S_{II}$ is not aperiodic by Proposition 
\ref{Ap*Gp}. It is known that the map $\GM_{J}: S \rightarrow \GM_{J}(S)$ is injective on subgroups of $J$ \cite[Chapter 4]{qtheory}. Let $J' = (J)\GM_{J}$. It follows easily that the image of $J \cap S_{II}$ is non-aperiodic and contained in $J' \cap (\GM_{J}(S))_{II}$. Therefore, 
$\GM_{J}(S)$ does not belong to $Ap*Gp$ by Theorem \ref{onepoint}. The second assertion now follows immediately from this and Theorem \ref{onepoint}. \Qed

$\GM$ inverse semigroups give important examples of $\GM$ semigroups that admit a 1-point flow. This is related to classical results in inverse semigroup theory. We review this for the reader's convenience and recast it in the language of flows.  For background in inverse semigroup theory, see \cite{Lawson, PetrichInv}.

Let $S$ be a $\GM$ inverse semigroup. Then its distinguished ideal with $I(S)$ is isomorphic to the Brandt semigroup $B_{n}(G) \approx M^{0}(G,\{1,\ldots n\}, \{1,\ldots n\},I_{n})$ for some $n$ and non-trivial group $G$ as these are the unique finite 0-simple semigroups. Let $B=\{1,\ldots n\}$. Then $S$ is a subsemigroup of $G \wr (B,SIM(B))$ and thus is isomorphic to a semigroup of $n \times n$ row and column monomial matrices over $G$. The congruence of the morphism from $S \rightarrow \RLM(S)$ using the inverse operation and \cite[Lemma 4.6.54]{qtheory} is the largest congruence on $S$ contained in $\mathcal{H}$. Therefore $\RLM(S)$ is the fundamental image of $S$ in the sense of inverse semigroup theory \cite{Lawson}. The decomposition $S$ divides $(G \wr Sym(B)) \times \RLM(S)$ then follows from a theorem of McAlister and Reilly \cite{McAlReilly} from 1977. In this case, the graph of the relational morphism from $S$ to $G \wr Sym(B)$ is an $E$-unitary cover of $S$ and the group $G \wr Sym(B)$ resolves the problem of the maximal group image of $S$ being the trivial group, since we assume $\GM$ semigroups have a zero. 

Let $(G \times B,S)$ is the corresponding $\GM$ transformation semigroup. Then $S$ acts by partial bijections on $G\times B$. It follows that the function from $f:\{1\} \rightarrow Rh_{B}(G)$ with $1f=(G \times B, \Pi)$, where $\Pi$ is the singleton partition of $G \times B$ is a 1-point flow.

As an example we look at the full wreath product $(G \times B, G \wr SIM(B))$. As mentioned above, we consider $G \wr SIM(B)$ to be the monoid of all $|B| \times |B|$ partial row and column monomial matrices over $G$. That is, the matrix of all $|B| \times |B|$ matrices with entries from $G\cup \{0\}$ that have at most 1 non-zero element in every row and column. The group of units of $G \wr SIM(B)$ is $G \wr Sym(B)$ the group of all monomial matrices over $G$. The semilattice of idempotents of $G \wr SIM(B))$ is isomorphic to $2^{B}$ under intersection, represented by all diagonal $|B| \times |B|$ matrices over $\{0,1\}$. Two matrices $M$ and $N$ are $\mathcal{R}$ ($\mathcal{L}$) related if and only if the set of non-zero columns (rows) of $M$ is equal to the set of non-zero columns (rows) of $N$. $M$ and $N$ are $\mathcal{J}$-related if and only if they have the same rank, which by definition is the number of non-zero rows. The maximal subgroup for the $\mathcal{J}$-class of rank $k$ is the group $G \wr Sym(k)$.

\bd

Let $M \in G \wr SIM(B)$. Then the domain, $Dom(M)$ of $M$ is the set of indices of its non-zero rows. The range $Ran(M)$ is the set of indices of its non-zero columns. 

\ed

Thus $M^{-1}M=I|_{Dom(M)}$ and $MM^{-1}=I|_{Ran(M)}$, where $I$ is the identity matrix.  We record this and more information in the following Lemma. The proofs are straightforward and are omitted.

\bl \label{GwrSIM}

Let $G$ be a non-trivial group and $B$ a set. Then the monoid of $(G \times B, G \wr SIM(B))$ is isomorphic to the monoid of all row and column monomial matrices over $G$. The following holds.

\begin{enumerate}
\item{$G \wr SIM(B)$ is a $\GM$ inverse monoid. The group of units is $G \wr Sym(B)$ the group of all $|B| \times |B|$ monomial matrices over $G$.}

\item{The semilattice of idempotents of $G \wr SIM(B)$ is the monoid of all diagonal matrices over $\{0,1\}$ and is isomorphic to $(2^{B},\bigcap)$}

\item{For $M,N \in G \wr SIM(B)$ we have:

\begin{enumerate}

\item{$M\mathcal{L}N$ if and only if $Dom(M)=Dom(N)$.}

\item{$M\mathcal{R}N$ if and only if $Ran(M)=Ran(N)$.}

\item{$M\mathcal{J}N$ if and only if the rank of $M$ is equal to the rank of $N$.}

\end{enumerate}}

\item{The maximal subgroup for the $\mathcal{J}$-class of rank $k$ is $G \wr Sym(k)$.}

\item{$\RLM(G \wr SIM(B))$ is isomorphic to $SIM(B)$ which is also the fundamental image of $G \wr SIM(B)$}

\item{There is a flow decomposition: $S \prec (G \wr Sym(B)) \times SIM(B)$.}

\item{The relational morphism $R:G \wr SIM(B) \rightarrow G \wr Sym(B)$ associated to this flow decomposition has graph 
$\{(M,N) \mid M \in G \wr SIM(B), N \in G \wr Sym(B), N|_{Dom{M}}=M\}$.}

\end{enumerate}
\el

For convenience we also write $R$ for the graph of the relational morphism in the last paragraph of Lemma \ref{GwrSIM}. $R$ is an $E$-unitary cover of $G \wr SIM(B)$. Let $\pi_{1}:R \rightarrow G \wr SIM(B), \pi_{2}:R \rightarrow G \wr Sym(B)$ be the projections. Then
$R = \pi_{1}^{-1}\pi_{2}$ is the canonical factorization of $R$ \cite{qtheory}. $\pi_{1}$ is an idempotent-separating morphism and it follows from functoriality, that the fundamental image of $R$ is also $SIM(B)$. $\pi_{2}$ is an idempotent-pure morphism, meaning that the semilattice of idempotents is a congruence class which by one of its definition proves that $R$ is an $E$-unitary semigroup. Note that from the point of view of complexity theory, idempotent-pure morphisms of inverse semigroups to groups are precisely the aperiodic morphisms to groups. 

It is well known that if $S$ is a finite inverse semigroup, the maximal group image of $S$ is isomorphic to the minimal ideal $K(S)$ of $S$. If $e$ is the unique idempotent in $K(S)$, then the map $\sigma:S \rightarrow K(S)$ such that $s\sigma = es$ is the map to the maximal group image of $S$.

\bl \label{Rstructure}

Let $R$ be the $E$-unitary cover of $G \wr SIM(B)$ described above. The following holds.

\begin{enumerate}

\item{The semilattice of idempotents of $R$ is isomorphic to $(2^{B},\bigcap)$}

\item{The morphism $\pi_{1}:R \rightarrow G \wr SIM(B)$ is idempotent separating. Therefore for any of the Green relations 
$\mathcal{K} \in \{\mathcal{R},\mathcal{L},\mathcal{H},\mathcal{J}\}$ we have $\faktor{R}{\mathcal{K}} \approx \faktor{G \wr SIM(B)}{\mathcal{K}}$. More precisely, $(M,N)\mathcal{K}(M',N')$ if and only if $M\mathcal{K} M'$.}

\item{ $SIM(B)$ is the fundamental image of $R$. }

\item{The maximal subgroup for the $\mathcal{J}$-class of rank $k$ is $G \wr Sym(k) \times (G \wr Sym(n-k))$. Here we regard $Sym(0)$ to be the empty set.}

\item{Consider the action of $(G \wr Sym(k))$ on the right of the $\mathcal{R}$-class of the idempotent $E_{k}=(1_{k},1)$ where $1_{k}$ is the diagonal matrix with a 1 precisely in rows $1 \ldots k$. Then the rank $k$ maximal subgroup $(G \wr Sym(k)) \times (G \wr Sym(n-k))$ is isomorphic to the right stabilizer of $E_k$.}

\end{enumerate}

\el

\proof An idempotent of $R$ is a pair $(E,1)$ where $E$ is an idempotent of $G \wr SIM(B)$. Therefore the projection to $E$ is an isomorphism between the semilattice of idempotents of $R$ to that of $G \wr SIM(B)$. Clearly, the semilattice of idempotents of 
$G \wr Sym(B)$ is the monoid of diagonal matrices over $\{0,1\}$ and is isomorphic to $(2^{B},\cap)$. This proves part 1. Since the projection from $E(R)$ to $E(G \wr Sym(B))$ is an isomorphism, we have proved that the projection from $R$ to $G \wr Sym(B)$ is idempotent separating. Since idempotent separating morphisms on inverse semigroups have congruences contained in $\mathcal{H}$, the claim about preservation of Green's relations is clear. This proves part 2. Furthermore, since the projection from $R$ to $G \wr SIM(B)$ is idempotent separating, it follows by functoriality that the fundamental image of $R$ is isomorphic to the fundamental image of $G \wr Sym(B)$, which is $SIM(B)$. This proves part 3.

As a representative of the rank $k$ $\mathcal{J}$-class we take $E_{k}=(1_{k},1)$ where $1_{k}$ is the diagonal matrix with 1 in rows $1, \ldots k$ 
and $0$ elsewhere. A pair $(M,N)$ is in the maximal subgroup $G_{k}$ with identity $E_{k}$ if and only if $MM^{-1}=M^{-1}M=1_{k}$. This is true if and only if $M$ is a monomial matrix of the form 
$\left[ \begin{array}{l r}
	A & 0 \\
	0 & 0
	
\end{array}\right]$ where $A \in G \wr Sym(k)$. Therefore $(M,N) \in G_{k}$ if and only if $M$ is as above and $N=\left[ \begin{array}{l r}
	A & 0 \\
	0 & B
	
\end{array}\right]$ where $B \in G \wr Sym(n-k)$. The map from such matrices to $(A,B)$ is an isomorphism from $G_{k}$ to $(G \wr Sym(k)) \times (G \wr Sym(n-k))$. This proves part 4. 

It is a straightforward calculation using part 4. to verify the statement in part 5. \qed

Now let $S$ be an arbitrary $\GM$ inverse semigroup. Then $I(S) \approx B_{n}(G)$ for a unique $n$ and a unique group $G$ up to isomorphism and $S$ is a subsemigroup of $G \wr SIM(B)$ containing $I(S)$ where $B =\{1, \ldots n\}$. We can restrict the relational morphism $R$ above to all pairs $(M,N)$ where $M \in S$. We give the appropriate generalization of Lemma \ref{Rstructure} to the general case. The proof is left to the reader.

\bl \label{GenRstructure}

Let $S$ be a $\GM$ inverse semigroup with $I(S) \approx B_{n}(G)$. Then $S$ is a subsemigroup of $G \wr SIM(B)$ containing $I(S)$.

Let $R_{S}=\{(M,N) \mid M \in S, N \in G \wr Sym(B), N|_{Dom(M)}=M\}$. Then $R_{S}$ is the graph of a relational morphism 
$S \rightarrow G \wr Sym(B)$ that is an $E$-unitary cover of $S$. The following holds.

\begin{enumerate}

\item{The semilattice of idempotents of $R_{S}$ is isomorphic to the semilattice of idempotents of $S$.}

\item\label{2}{The projection $\pi_{1}:R_{S} \rightarrow S$ is a surjective idempotent separating morphism. Therefore for any of the Green relations 
$\mathcal{K} \in \{\mathcal{R},\mathcal{L},\mathcal{H},\mathcal{J}\}$ we have $\faktor{R_{S}}{\mathcal{K}} \approx \faktor{S}{\mathcal{K}}$. More precisely, $(M,N)\mathcal{K}(M',N')$ if and only if $M\mathcal{K} M'$.}

\item{$RLM(S)$ is the fundamental image of $R_{S}$. The projection of $R_{S}$ to $G \wr Sym(B)$ has image a subgroup $H_{S}$ of $G \wr Sym(B)$. $R_{S}$ is a subsemigroup of $H_{S} \times RLM(S)$.}

\item{Let $J$ be a $\mathcal{J}$-class of $S$ and $J'$ be the $\mathcal{J}$-class of $R_{S}$ corresponding to $J$ from part \ref{2}. The maximal subgroup for $J'$ is the subgroup of $H_{S}$ that stabilizes an idempotent of $J'$.}

\end{enumerate}

\el

The results in this section on inverse semigroups above are a classical part of inverse semigroup theory. An inverse monoid $M$ is called
{\em factorizable} if $M=EU$, where $E$ is the semilattice of idempotents of $M$ and $U$ is the group of units of $M$. This means that every element of $M$ is less than or equal to some element in $U$ in the natural partial order of $M$. If $M$ is factorizable, then 
$R_{M}=\{(m,g)\mid m \in M, g \in U, m \leq g\}$ is the graph of a relational morphism $M \rightarrow U$ that is an $E$-unitary cover of $M$. Thus each factorizable inverse monoid $M$ naturally defines an $E$-unitary cover over its group of units. Conversely, define a strict embedding of an inverse semigroup $S$ into a factorizable inverse monoid $M$ to be an embedding $\iota:S \rightarrow M$ such that for each 
$g \in U$, where $U$ is the group of units of $M$, there is $s \in U$ such that $s\iota \leq g$. Then the relational morphism $S \rightarrow U$ that sends $s \in S$ to the set of $\{g \in G|s \leq g\}$ has graph an $E$-unitary cover of $M$ with maximal group image $U$. This is precisely the construction of $R_{S}$ in Lemma \ref{GenRstructure} since $G \wr SIM(B)$ is a factorizable inverse monoid with group of units $G \wr Sym(B)$. For the precise connection between $E$-unitary covers of an inverse semigroup $S$ with maximal group image $U$ and strict embeddings into factorizable inverse monoids with group of units $U$ see Section 8.2 of \cite{Lawson}.

We have seen that 1-point flows are intimately related to the pseudovariety $Ap*Gp$. We introduce the notion of a set-trivial flow and relate it to the pseudovariety $Gp*Ap$.

\bd

Let $(G\times B,S)$ be a $\GM$ transformation semigroup. Let $(Q,T)$ be an aperiodic transformation semigroup. A set-trivial flow is a flow $Q \rightarrow Rh_{G}(B)$ such that for each $q \in Q$, the set of the SPC $qF$ is a singleton.

\ed

\bt
Let $(G \times B,S)$ be a $\GM$ transformation semigroup. The following conditions are equivalent.

\begin{enumerate}

\item{$S \in Gp*Ap$}

\item{$\RLM(S)$ is aperiodic.}

\item{$(G \times B,S)$ admits a set-trivial flow.}
\end{enumerate}

\et

\proof The equivalence of item 1. and item 2. follows from \cite{complex1/2}. Assume now that $\RLM(S)$ is aperiodic. Then $(B,\RLM(S))$ is an aperiodic transformation semigroup. The map $f:B \rightarrow Rh_{B}(G)$ with $bf=b/\langle 1 \rangle $ is clearly a set-trivial flow. This proves that item 2. implies item 3.

Now assume that $(G \times B,S)$ admits a set-trivial flow over the aperiodic transformation semigroup $(Q,T)$. This means that there is a map
from $f:Q \rightarrow Rh_{G}(B)$ with $qf =b/\langle 1 \rangle $ for some $b \in B$ and such that every $b \in B$ is the image of some $q \in Q$. Therefore for each $q \in Q, s \in S$, we have that if $(qf)s$ is defined, then  $(qf)s = g_{q,qs}(qs)f$ for a unique $g_{q,qs} \in G$. Assigning $s$ to the $Q \times Q$ row-monomial matrix whose entries are the non-zero $g_{q,qs}$ then shows that $(G \times B,S)$ divides $G \wr (Q,T)$ and thus $S \in Gp*Ap$. This completes the proof. \Qed

\section{Slices, Flows and the Presentation Lemma: A Unified Approach} \label{unified}


There are three versions of the Presentation Lemma and its relation to Flow Theory in the literature \cite{AHNR.1995}, \cite[Section 4.14]{qtheory}, \cite{Trans}. These have very different formalizations and it is not clear how to pass from one version to another. The terminology is different as well. For example, the definition of cross-section in each of these references, as well as the one we use in this paper is different.  The purpose of this section is the following Theorem that gives a unified approach to these topics. It summarizes known results in the literature and is meant to emphasize the strong connections between slices, flows and the presentation lemma and the corresponding direct product decomposition. We use some of the results that we've proved in previous sections. For background on the derived transformation semigroup and the derived semigroup theorem see \cite{Eilenberg, qtheory}.

\bt \label{uniform}

Let $(G \times B,S)$ be $\GM$ and assume that $\RLM(S)c \leq n$. Then the following are equivalent:

\begin{enumerate}

\item{$Sc \leq n$.}

\item{There is an aperiodic relational morphism $\Theta:S \rightarrow H\wr T$, where $H$ is a group and $Tc \leq n-1$.}

\item{There is a relational morphism $\Phi:S \rightarrow T$ where $Tc \leq n-1$ and such that the Derived Transformation semigroup $D(\Phi)$ is in $Ap*Gp$.}

\item{There is a relational morphism $\Phi:S \rightarrow T$ where $Tc \leq n-1$ such that the Tilson congruence $\tau$ on the Derived Transformation semigroup $D(\Phi)$ is a cross-section.}

\item{$(G \times B,S)$ admits a flow from a transformation semigroup $(Q,T)$ with $(Q,T)c \leq n-1$.}

\item{$S \prec (G \wr Sym(B) \wr T) \times \RLM(S)$ for some transformation semigroup $T$ with $Tc \leq n-1$.}

\end{enumerate}

\et

\proof  Assume that $Sc \leq n$. Then $\RLM(S)c \leq n$ and $S \prec A \wr H \wr T$ where $A$ is aperiodic, $H$ is a group and $Tc \leq n-1$. Therefore, there is a division $d:S \rightarrow A \wr H \wr T$. For 2. we take $\Theta = d\pi$ where $\pi: A \wr H \wr T \rightarrow H \wr T$ is the projection. and 3. follows by taking $\Phi=\Theta\pi'$ where $\pi'$ is the projection from $H \wr T$ to $T$ and the Derived Semigroup Theorem. The Derived Semigroup Theorem also imply that 2. and 3. imply 1. Therefore 1., 2. and 3. are equivalent.

By Theorem \ref{onepoint} it follows that 3. is equivalent to 4. The equivalence of 4. and 5. follows from Section 3 of \cite{Trans} and section 4.14 of \cite{qtheory}. That 5. implies 6. follows from the Slice Theorem and Theorem \ref{GwrSymB}.  Since the complexity of a direct product is the maximum of its factors, 6. implies 1. by our assumption on $\RLM(S)c$. This completes the proof. \Qed

Since the emphasis in this paper is on aperiodic flows, we restate Theorem \ref{uniform} in the case of aperiodic flows.

\bt

Let $(G \times B,S)$ be $\GM$ and assume that $\RLM(S)c \leq 1$. Then the following are equivalent:

\begin{enumerate}

\item{$Sc=1$.}

\item{There is an aperiodic relational morphism $\Theta:S \rightarrow H\wr T$, where $H$ is a group and $T$ is aperiodic.}

\item{There is a relational morphism $\Phi:S \rightarrow T$ where $T$ is aperiodic such that the Derived transformation semigroup $D(\Phi)$ is in $Ap*Gp$.}

\item{There is a relational morphism $\Phi:S \rightarrow T$ where $T$ is aperiodic such that the Tilson congruence $\tau$ on the Derived transformation semigroup $D(\Phi)$ is a cross-section.}

\item{$(G \times B,S)$ admits an aperiodic flow.}

\item{$S \prec (G \wr Sym(B) \wr T) \times \RLM(S)$ for some aperiodic semigroup $T$.}

\end{enumerate}

\et

\section{Flows for Some Small Monoids}\label{smallmon}

In the following two sections we will look at semigroups that have flows over semigroups that are monoids consisting of a minimal ideal that is a right-zero semigroup with an identity adjoined. The right-zero semigroup with $k$ elements is denoted $RZ(k)$ and the monoid obtained from this semigroup by adjoining an identity element is denoted by $RZ(k)^{1}$. Thus resolution semigroups will be subsemigroups of $H \wr RZ(k)^{1}$ for some group $H$ and integer $k>0$. We first elucidate the structure of these semigroups. Recall that a completely regular semigroup is a semigroup that is a union of its subgroups. A small monoid is a monoid $M$ consisting of a group of units and a $0$-minimal ideal $I(M)$.

\bt \label{gpwreathrz}

\begin{enumerate}

\item{A finite semigroup $S$ is a $\RLM$ simple semigroup if and only if $S$ is isomorphic to $RZ(k)$ for some $k>0$.}

\item{A finite monoid $M$ is an aperiodic completely regular small $\RLM$ monoid if and only if $M$ is isomorphic to $RZ(k)^{1}$ for some $k$.}

\item{Let $H$ be a non-trivial group. Let $S(H,k)$ be the semigroup of the transformation semigroup $H \wr (\{1,\ldots k\},RZ(k))$. Then $S(H,k)$ is the largest $\GM$ simple semigroup with $k$ 
$\mathcal{L}$-classes and maximal subgroup $H$. That is, every simple semigroup with these properties is a subsemigroup of $S(H,k)$.}

\item{$S(H,k)$ has a Rees matrix representation $S(H,k) \approx M(H,\{f:\{1,\ldots k\}\rightarrow H \mid f(1)=1_{H}\},\{1,\ldots k\},C)$ where $C(i,f)=if$. If we take the first column of $C$ to be the constant function to $1_H$, then $C$ is the biggest matrix with entries in $H$, $k$ rows and whose first row and first column are all the identity of $H$ and such that no two rows nor columns of $C$ are proportional. $C$ is the matrix with $k$ rows and the maximal number of columns with these properties.}

\item{Let $H$ be a non-trivial group. The monoid $M(H,k)$ of $H \wr (\{1,\ldots k\},RZ(k)^{1})$ is a small $\GM$ completely regular monoid with minimal ideal $S(H,k)$ and group of units $H^{k}$.}

\item{Let $M(H,k)=H^{k} \cup S(H,k)$ as in the previous item. Then the right Sch\"{u}tzenberger representation of $H^{k}$ on $S(H,k)$ sends
 $(h_{1},\ldots, h_{k}) \in H^{k}$ to the diagonal $k \times k$ matrix with $i^{th}$ entry $h_{i}$.}

\item{Let $M(H,k)=H^{k} \cup S(H,k)$. Then the left Sch\"{u}tzenberger representation of $H^{k}$ on $S(H,k)$ sends
 $(h_{1},\ldots, h_{k}) \in H^{k}$ to the $|H|^{k-1} \times |H|^{k-1}$ matrix all of whose non-zero entries are $h_{1}$ and whose entry in column $i$ is the unique column $C_{j}$ such that $C_{j}h_{1}=C_{i}$.}

\end{enumerate}

\et

\proof

Items 1. and 2. are immediate from the definitions. We turn to items 3. and 4. As a set we can let 
$S(H,k)=H^{k} \times RZ(\{1,\ldots k\})$, where we write $\bar{i}$ for the constant function on $\{1,\ldots k\}$ with value $i$. With these identifications multiplication in $S(H,k)$ is given by: 
$$((h_{1},\ldots, h_{k}),\bar{i})((h_{1}^{'},\ldots, h_{k}^{'}),\bar{j})=((h_{1}h_{i}',\ldots, h_{k}h_{i}^{'}),\bar{j})$$.

It follows that $((h_{1},\ldots, h_{k}),\bar{i})^{n}=((h_{1},\ldots, h_{k}),\bar{i})$, where $n$ is the order of $h_{i}$ and thus $S(H,k)$ is a completely regular semigroup. Furthermore,  
$((h_{1},\ldots, h_{k}),\bar{i})=((h_{1}h_{i}'^{-1},\ldots, h_{k}h_{i}'^{-1}),\bar{i})((h_{1}^{'},\ldots, h_{k}^{'}),\bar{i})$ and thus for 
every $1 \leq i \leq k$ it follows that $L_{i}=H^{k} \times \{i\}$ is an $\mathcal{L}$-class. Since $L_{i}L_{j}=L_{j}$, we have that $\mathcal{L}$ is a congruence on $S(H,k)$ with quotient $RZ(k)$. Therefore, $S(H,k)$ is a completely simple semigroup. 

Note that 
$((h_{1},\ldots, h_{k}),\bar{i})$ is an idempotent of $S(H,k)$ if and only if $h_{i}=1_{H}$. It follows that $S(H,k)$ has $|H|^{k-1}$ $\mathcal{R}$-classes. Consider the idempotent $e=((1_{H},\ldots , 1_{H}),\bar{1})$. Since $S(H,k)$ is completely simple, the maximal subgroup $G_{e}$ is equal to $eS(H,k)e$. A straightforward calculation shows that $G_{e}=\{((h,h,\ldots ,h),\bar{1}) \mid h \in H\} \approx H$. Let
$A=\{((1_{H},h_{2},\ldots,h_{k}),\bar{1}) \mid (h_{2},\ldots,h_{k})\in H^{k-1}\}$ and $B=\{((1_{H},\ldots , 1_{H}),\bar{i}) \mid i=1,\ldots k\}$. Then $A$ ($B)$ is both the set of idempotents in the $\mathcal{L}$-class ($\mathcal{R}$-class) of $e$ and a set of representatives for the $\mathcal{R}$-classes ($\mathcal{L}$-classes) of $S(H,k)$. The proof of the Rees Theorem then gives us the Rees representation
$S(H,k) \approx M(H,A,B,C)$ where $C$ is the $B \times A$ matrix whose entry in place
$((1_{H},\ldots , 1_{H}),\bar{i}), ((1_{H},h_{2},\ldots,h_{k}),\bar{1})$ is 

$((1_{H},\ldots , 1_{H}),\bar{i})((1_{H},h_{2},\ldots,h_{k}),\bar{1}) = ((h_{i},h_{i},\ldots , h_{i}),\bar{1})$. We identify $A$ with the set
of functions $f:\{1,\ldots , k\} \rightarrow H$ with $1f =1_{H}$ by sending such an $f$ to $((1_{H},2f,\ldots , kf),\bar{1})$ and $B$ with the set $\{1,\ldots , k\}$ in the obvious way and $G_{e}$ with $H$ by sending $(h,h,\ldots , h),\bar{1})$ to $h$. With these identifications, the Rees matrix representation we've found is exactly as stated in item 4. Finally any $\GM$ simple semigroup $S$ with $k$ 
$\mathcal{L}$-classes and maximal subgroup $H$ is known to embed into $H \wr \RLM(S)$. It is clear that for such an $S$, $\RLM(S) \approx RZ(k)$ and thus $S$ embeds into $S(H,k)$. This completes the proof of item 3.

The monoid $M(H,k)$ of $(H \wr RZ(k)^{1})$ is the union of $S(H,k)$ and the set $\{((h_{1},\ldots , h_{k}),1) \mid (h_{1},\ldots , h_{k})\}$. Clearly this latter set is the group of units of $M(H,k)$ and is isomorphic to $H^{k}$. By item 4., it follows that $S(H,k)$ is the minimal ideal. We have proved item 5. The proofs of items 6. and 7. are straightforward computations. This completes the proof of the Theorem.
\Qed

\be

Let $H = Z_2= \{1,-1\}$. Then $S(Z_{2},2)$ is the simple semigroup $M(Z_{2},\{1,2\},\{1,2\},C)$ where $C=
\begin{bmatrix}

1 & 1 \\
1 & -1

\end{bmatrix}
$. $M(Z_{2},2) = S(Z_{2},2) \cup (Z_{2} \times Z_{2})$, where $(x,y)\in Z_{2} \times Z_{2}$ acts in the right Sch\"{u}tzenberger representation on $S(Z_{2},2)$ by the diagonal matrix $\begin{bmatrix}

x & 0 \\
0 & y

\end{bmatrix}$. It is known that the ten element monoid, that is the reverse of $S(H,2) \cup \{(1,-1)\}$ is the smallest semigroup of complexity 2. $M(H,k)$ has complexity 1 for all $H$ and $k$ as it is explicitly defined as the wreath product of an aperiodic semigroup and a group.

\ee

Small monoids play an important part in relations between semigroups and combinatorics \cite{Lallement, Pollatchek} and in representation theory of finite monoids \cite{Putcharep3}. Tilson showed that complexity is decidable for small monoids and more generally for semigroups with at most 2 non-zero  $\mathcal{J}$-classes \cite{2J}. Historically this paper was the first that used techniques related to the presentation lemma and flows. Before giving the details of the theorem we give some examples. We begin with a useful Lemma. 

\bl\label{smalltypeII}
Let $S$ be a small $\GM$ monoid with group of units $H$ and 0-minimal ideal $M^{0}(G,A,B,C)$. Then the following conditions are equivalent:

\begin{enumerate}

\item{$S$ divides $T \wr K$ for some aperiodic semigroup $T$ and group $K$.}

\item{$S_{II}$ is aperiodic.}

\item{The idempotent generated subsemigroup $IG(S)$ is aperiodic.}

\item{The structure matrix $C$ is normalizable to $\{0,1\}$.}

\item{$S$ admits a $1$-point flow.}

\item{$S$ divides $G \wr Sym(B) \times \RLM(S)$.}

\end{enumerate}
\el

\proof
It is known that item 1. is equivalent to item 2. and that item 3. is equivalent to item 4. \cite{Graham}, \cite[Section 4.13]{qtheory}. Clearly item 2. implies item 3. We prove that item 4. implies item 1. 

If all of the entries of $C$ belong to $\{0,1\}$, then $M^{0}(1,A,B,C)$ is a subsemigroup of $M^{0}(G,A,B,C)$. Furthermore, the map sending $(g,(a,1,b))$ to $(a,g,b)$ is a surjective morphism from $G\times M^{0}(1,A,B,C)$ to $M^{0}(G,A,B,C)$. It follows that the submonoid 
$N =\{1\} \cup M^{0}(G,A,B,C)$ divides $G \times (M^{0}(1,A,B,C) \cup \{1\}$. Therefore $N$ divides the direct product and hence the wreath product of an aperiodic semigroup and a group. The group of units $H$ acts by conjugation on $N$ and the map from the corresponding semidirect product to $M$ sending $(1,h)$ to $h$ and $((a,g,b),h)$ to $(a,g,b)h$ is a surjective morphism. Since semidirect products embed into wreath products, item 1. is proved.

We have proved that items 1. through 4. are equivalent. Theorem \ref{onepoint} immediately implies that item 2. is equivalent to item 5. Finally, the equivalence of item 5. and item 6. follows from Theorem \ref{onepoint}. This completes the proof of the Lemma. \qed

\brm
It is not hard to see more generally that if $S$ is a small monoid, then the maximal subgroups of $S_{II}$ are the trivial subgroup of $H$ and the normal closures in $G$ of the maximal subgroups of $IG(S)$. See \cite{lowerbounds2, Graham}. We do not need this result in the paper.
\erm

\brm 

Example 5.6 of \cite{Mcalisterschutz} gives an example of a regular semigroup with 2 non-zero $\mathcal{J}$-classes in which $IG(S)$ is aperiodic, but $S_{II}$ is not aperiodic. Thus Lemma \ref{smalltypeII} can not be generalized beyond small monoids, even for semigroups with 2 non-zero $\mathcal{J}$-classes.

\erm
 



We now give some examples of small monoids. One will have a $1$-point flow, one has complexity 1, but admits no 1-point flow and one has complexity 2.

Let $S$ be a $\GM$ semigroup with 0-minimal ideal $M^{0}(G,A,B,C)$. As $S$ embeds into $G \wr \RLM(S)$, $S$ has a faithful representation by 
$|B| \times |B|$ row monomial matrices over $G$. It is convenient to define such matrices as partial functions on $B$ labelled by elements of $G$. For example if the row monomial matrix has the element $g$ in position $(i,j)$ we write $i \rightarrow gj$, if $g \neq 1$ and 
$i\rightarrow j$ if $g = 1$.  

\be

We give an example of a non-inverse semigroup with a 1-point flow. Let $S$ be the $0$-simple semigroup 
$M^{0}(G,A,B,C)$ where $G=Z_{2}, A = \{a_{1},a_{2},a_{3},a_{4}\}, B=\{b_{1},b_{2},b_{3},b_{4}\}$ and 
$$C = \left[\begin{array}{c c c c} 1&1&0&0 \\ 0&1&1&0 \\ 0&0&1&1\\ 1&0&0&1\end{array}\right]$$. 

Let $Z_4$ be the cyclic group of order 4 generated by an element $z$. For $k \in \{0,1,2,3\}$, define $(a_{i},g,b_{j})z^{k}=(a_{i},g,b_{j+k \operatorname{mod} 4}, z^{k}(a_{i},g,b_{j})=
(a_{i-k \operatorname{mod} 4},g,b_{j})$. Then with this definition as multiplication of $Z_4$ on $S$, $M = Z_{4} \cup S$ is a $\GM$ small monoid with group of units $Z_4$ and with $I(M)=S$. Therefor $M$ divides $(G \wr Sym(4)) \times \RLM(S)$ and since $\RLM(S)c=1$, we have $Mc=1$.


\ee

Before giving more examples, we state Tilson's Theorem on the complexity of small monoids. In fact, the theorem give a decidability criterion for all semigroups with 2 non-zero $\mathcal{J}$-classes. The general case can be reduced to that of small monoids. See \cite{2J}, 
\cite[Section 4.15]{qtheory} for details. From the depth decomposition theorem \cite{TilsonXI}, a small monoid has complexity at most 2. Thus we need a criterion to decide whether the complexity is 1 or not.

Let $S = H \cup M^{0}(G,A,B,C)$ be a small monoid. We identify $B$ with the set of $\mathcal{L}$-classes of $M^{0}(G,A,B,C)$. Then $H$ acts on the right of $B$ by $L\cdot h = Lh$, where $L$ is an $\mathcal{L}$-class of $M^{0}(G,A,B,C)$. This is the restriction of the action defining $\RLM(S)$ to $H$. Then for each such $L$, $H \cup LH \cup \{0\}$ is a submonoid of $M$ called a right-orbit monoid. The collection of all $LH$ are the orbits of the action of $H$ on $B$. Let $k=k(S)$ be the number of such orbits. Let $B_{1}, \ldots B_{k}$ be the orbits of $H$ on $B$. Then $G \times B_{1}, \ldots , G \times B_{k}$ is a partition of $G \times B$. Furthermore the right-orbit monoids are exactly the monoids $H \cup A \times G \times B_{i}, i=1, \ldots k$. Notice that $A \times G \times B_{i} = M^{0}(G,A,B_{i},C_{i})$, where $C_{i}$ is the restriction of $C$ to $B_{i} \times A$. This is a not necessarily regular Rees matrix semigroup.

\bt\label{2J}

Let $S = H \cup M^{0}(G,A,B,C)$ be a small $\GM$ monoid. Then $Sc = 1$ if and only if each right-orbit monoid $H \cup LH \cup \{0\}$ has aperiodic idempotent generated subsemigroup. 

\et

The proof of Theorem \ref{2J} is stated using the Presentation Lemma \cite{AHNR.1995}. We recast it in the language of flows. The proof is given in \cite{2J} and in a slightly different form in Section 4.15 of \cite{qtheory} in the language of the Presentation Lemma. To complete the proof one uses the equivalence between the Presentation Lemma and Flow Theory given in Section 3 of \cite{Trans}.

\bt\label{2Jflow}

Let $S$ be a small monoid with $k=k(S)$ left orbits. Assume that every right-orbit monoid $H \cup LH \cup \{0\}$ has aperiodic idempotent generated subsemigroup. Let $RZ(k)^{1}$ be the right-zero semigroup $RZ(k)$ with an identity adjoined. We have the transformation monoid
$(\{1 \ldots k\},RZ(k)^{1})$. Then the following defines a flow $\Phi:\{1 \ldots k\} \rightarrow Rh_{B}(G)$.

Let $i\Phi = (G \times B_{i},\Pi_{i})$, where $\Pi_{i}$ is the partition defined by the equivalence relation $(g,b)\Pi_{i}(g',b')$ if and only if there is an element $s$ in the idempotent generated subsemigroup of $A \times G \times B_{i}$ with $(g,b)s = (g',b')$.
Therefore $S$ divides $(G \wr Sym(B) \wr RZ(k)^{1}) \times \RLM(S)$.

\et

\be

We give an example of a small monoid of complexity 1, that has no 1-point flow. Let $S$ be the $0$-simple semigroup 
$M^{0}(G,A,B,C)$ where $Z_{2}=\{1,-1\}, A = \{a_{1},a_{2},a_{3},a_{4}\}, B=\{b_{1},b_{2},b_{3},b_{4}\}$ and 
$$C = \left[\begin{array}{c c c c} 1&1&0&0 \\ 0&1&1&0 \\ 0&0&1&1\\ 1&0&0&-1\end{array}\right]$$. Then it is known that $C$ is not-normalizable to $\{0,1\}$ and thus $IG(S)$ is not aperiodic. More concretely, $(1,1,1)(2,1,2)(3,1,3)(4,-1,4)(1,1,1) = (1,-1,1)$ is a non-trivial group element in $IG(S)$. Therefore no $\GM$ semigroup with $S$ as the 0-minimal ideal has a 1-point flow.

Let $s$ be defined by $b_{1} \rightarrow -b_{3}, b_{2}\rightarrow b_{4}, b_{3} \rightarrow b_{1}, b_{4} \rightarrow -b_{2}$. The subgroup of $Z_{4} \wr Sym(4)$ generated by $s$ is a cyclic group $Z_{4}$ of order 4. Direct calculation shows that $M = Z_{4} \cup S$ is a small monoid with group of units $Z_{4}$ and $I(M) \approx S$.

There are two right-orbit monoids with orbits $A \times Z_{4} \times \{b_{1},b_{3}\}$ and $A \times Z_{4} \times \{b_{2},b_{4}\}$. Again, direct calculation proves that each of these have aperiodic idempotent generated semigroups. Therefore, by Theorem \ref{2J}, we have $Mc=1$. By 
Theorem \ref{2Jflow}, $M$ divides $(Z_{4} \wr Sym(4) \wr RZ(2)^{1}) \times RLM(M)$.
\ee

\be

Let $S$ be the 0-simple semigroup in the previous example. Let $x$ be defined by $b_{1} \rightarrow b_{4}, b_{2} \rightarrow b_{1}, b_{3} \rightarrow b_{2}, b_{4} \rightarrow -b_{3}$. Then $x$ generates a subgroup of $Z_{2} \wr Sym(4)$ of order 8. Direct calculation shows that $M = Z_{8} \cup S$ is a small monoid with group of units $Z_{8}$ and $I(M) \approx S$. There is only 1 orbit of right-orbit monoids and since the idempotent generated subsemigroup of $S$ is not aperiodic, we have $Mc = 2$ by Theorem \ref{2J}.

\ee

\section{$\GM$ semigroups built from character tables of Abelian Groups}\label{chartab}

In this section we build a $\GM$ semigroup $S(n)$ of complexity 1 using the character theory of a cyclic group of order $2^n$. We show that $S(n)$ has a flow over $RZ(n)^1$. 

We write the cyclic group of order $n$ with generator $x$ multiplicatively. That is, $Z_{n}=\{x^{i} \mid i = 0,\ldots , n-1\}$.
 We define the ``abstract'' character table of $Z_{n}$ to be the matrix $C_n$ with rows and columns indexed by $\{0,\ldots , n-1\}$ and such that 
$C_{n}(k,l)=x^{kl}$. If we substitute the complex $n^{th}$ primitive root of unity, $e^{2\pi i/n}$ for $x$ 
 we obtain the complex character table $\chi_{n}$ of $Z_{n}$. Thus in $C_n$ we are effectively identifying $Z_{n}$ with its complex character group. 

We define $Ch_{n}$ to be the simple semigroup $M(Z_{n},\{0,\ldots , n-1\},\{0,\ldots , n-1\},C_{n})$. Since the rows and columns of $Ch_n$
are not proportional, $Ch_{n}$ is a $\GM$ simple semigroup with maximal subgroup $Z_n$. By Theorem \ref{gpwreathrz}, $Ch_n$ is a subsemigroup of the semigroup $S(Z_{n},n)$ of the wreath product $Z_{n} \wr (n,RZ(n))$. 

Let $X$ be the $n \times n$ permutation matrix with indices 
$\{0,\ldots , n-1\}$ such that $X(i,i+1)=1$, where we take $i$ modulo $n$. That is, $X$ is the permutation matrix of the cycle $(01\ldots n-1)$ acting on the right. Let $Y$ be the $n \times n$ diagonal monomial matrix with $Y(i,i) = x^{i}$. Direct matrix multiplication shows that
$XC_{n} = C_{n}Y$. Effectively, shifting the rows up 1 is the same as multiplying column $i$ by $x^{i}$. Therefore, the pair $(X,Y)$ is in the translational hull of $Ch_{n}$ and we can form the monoid $R_{n} = Ch_{n} \cup Z_{n}$ where the generator $x$ of $Z_{n}$ acts via $X$ on the right of $Ch_{n}$ and by $Y$ on the left. The monoid $Res_{n}=R_{n}^{op}$ is then a submonoid
of $M(Z_{n},n)$ by Theorem \ref{gpwreathrz}. We need the opposite monoid because we want a diagonal left action of $x$ on $S(Z_{n},n)$ since SPCs are proportional on the left. We will see below that $Res_{n}$ is the resolution semigroup of the monoid $S_{n}$ that we now define.

We remark that there is a natural representation theoretic meaning to this construction. If we replace $C_{n}$ by the complex character table 
$\chi_{n}$ as explained above, then we can rewrite the equation $XC_{n}=C_{n}Y$ as $Y = C_{n}^{-1}XC_{n}$, where we also replace $x \in Y$ by $e^{2\pi i/n}$. This reflects the well-known fact that the eigenvectors of $\chi_{n}$ give a basis that diagonalizes the right regular representation of $Z_n$ as a direct sum of one copy of each irreducible representation.

We now define a semigroup $S_{n}$ for each $n>0$ with the following properties:

\begin{enumerate}

\item{$S_{n}$ is a $\GM$ semigroup whose 0-minimal $\mathcal{J}$-class has maximal subgroup $Z_{2^{n}}$.}

\item{$S_{n}$ has a flow over the semigroup $RZ(n)^{1}$.}

\item{The semigroup $Res_{n}$ defined above is the resolution semigroup for this flow.}

\end{enumerate} 

The semigroup $S(1)$ is discussed in Example 5.11 of \cite{deg2part2}. For ease of notation and to look at a concrete example we define $S(2)$. We give the easy generalization to $S(n)$ for all $n$ after this.

\be \label{S2}

Let $\Gamma_{8}$ be the incidence matrix of an 8-cycle. More precisely, $\Gamma_{8}$ is the $8 \times 8$ matrix whose rows are labeled by the vertices $\{1,\ldots 8\}$ and columns labeled by the 8 edges, which for convenience we also use $\{1,\ldots 8\}$. Therefore the $i^{th}$ column of $\Gamma_{8}$ has a 1 in rows $i$ and $i+1$ and 0 elsewhere. Let $M_{8}$ be the $8 \times 16$ matrix who in block form is the matrix
$\begin{bmatrix} \Gamma_{8} & | & I_{8} \end{bmatrix}$ where $I_{8}$ is the identity matrix of size 8.

By Lemma 4.1 of \cite{cremona}, the translational hull of the 0-simple semigroup over the trivial semigroup with structure matrix $M_{8}$ consists of all partial functions $f$ on the vertices of the 8-cycle such that the inverse image of a vertex or edge in the image of $f$ is also either a vertex or an edge. This makes it easy to check if such a function belongs to the translational hull of this 0-simple semigroup. More generally, if we let $M_{8}(G)$ be the 0-simple semigroup with maximal subgroup $G$ and structure matrix $M_{8}$, Lemma 4.1 of \cite{cremona} gives a necessary and sufficient condition for a row monomial matrix over $G$ to be in the translational hull of $M_{8}(G)$. We use this criterion without mentioning it in the definition of $S(2)$ and later $S(n)$.

We now define $S(2)$ as a subsemigroup of the translational hull of $M_{8}(Z_{4})$. The generators are all the elements of the 0-minimal ideal $M_{8}(Z_{4})$ together with the following elements:

\begin{enumerate}

\item{Let $a$ be the permutation $(1357)(2468)$ corresponding to adding 2 modulo 8. By our convention, this is shorthand for the row-monomial permutation matrix corresponding to $a$.}

\item{Define $b$ = $8 \rightarrow 1$, $7 \rightarrow 2$. Notice that $a$ maps an element in its domain from an orbit of $a$ to an element in the other orbit.}

\item{We now use the abstract character table $C_4 = \begin{bmatrix} 1 & 1 & 1 & 1 \\ 1&x&x^{2}&x^{3} \\ 1&x^{2}&1&x^{2} \\ 1&x^{3}&x^{2}&x \end{bmatrix}$ to define 3 more generators:
\begin{enumerate}

\item{$s_{x}= 1 \rightarrow 2, 2 \rightarrow x4, 3 \rightarrow x^{2}6, 4 \rightarrow x^{3}8$. That is, $s_{x}$ weights the function that sends 
$1\rightarrow 2, 2\rightarrow 4, 3\rightarrow 6, 4\rightarrow 8$ by the character corresponding to $x$.}

\item{$s_{x^{2}}$ weights $1\rightarrow 2, 2\rightarrow 4, 3\rightarrow 6, 4\rightarrow 8$ by the character corresponding to $x^2$. Thus, 
$s_{x^{2}}= 1 \rightarrow 2, 2 \rightarrow x^{2}4, 3 \rightarrow 6, 4 \rightarrow x^{2}8$.}

\item{$s_{x^{3}}$ weights $1\rightarrow 2, 2\rightarrow 4, 3\rightarrow 6, 4\rightarrow 8$ by the character corresponding to $x^3  $. Thus, 
$s_{x^{2}}= 1 \rightarrow 2, 2 \rightarrow x^{3}4, 3 \rightarrow x^{2}6, 4 \rightarrow x8$.}

\end{enumerate}}

\item{We now define $S_2$ to be the semigroup generated by the 0-simple semigroup $M_{8}(Z_4)$, where we take $x$ as the generator of the distinguished copy of $Z_4$ in this semigroup and the elements $a,b,s_{x},s_{x^{2}},s_{x^{3}}$ defined above.}

\end{enumerate}

We first gather some facts about $S_{2}$.

\begin{enumerate}

\item{$S_{2}$ is a subsemigroup of the translational hull of $M_{8}(Z_{4})$. As mentioned above, this is easily justified by using Lemma 4.1 of \cite{cremona} and Lemma 4.1 of \cite{cremona}. It follows that $S_{2}$ is a $\GM$ semigroup. Furthermore, 
$I(S_{2}) \approx M^{0}(Z_{4},\{1,\ldots 16\}, \{1 \ldots 8\}, M_{8})$.}

\item{$S_{2}$ does not have a one-point flow. By Theorem \ref{onepoint} we must show that $I(S_{2}) \cap (S_{2})_{II}$ is not aperiodic, which we now do. Firstly, by Graham's Theorem \cite{Graham} it follows that $IG(S_{2}) \cap I(S_{2}) \approx M^{0}(\{1\},\{1,\ldots 16\}, \{1 \ldots 8\}, M_{8})$ is aperiodic, where $IG(S_{2})$ is the idempotent generated subsemigroup of $S_{2}$. 

Now note that $(4,x^{3},2)s_{x}(4,x^{3},2)=(4,x^{3},2)$. Therefore $(4,x^{3},2)ts_{x} \in (S_{2})_{II}$ for any $t \in (S_{2})_{II}$. Let
$t=(1,1,1) = t^{2} \in  (S_{2})_{II}$. Thus $(4,x^{3},2)ts_{x}=(4,x^{3},2)(1,1,1)s_{x}=(4,x^{3},2) \in (S_{2})_{II}$. By the above, 
$(2,1,4) \in IG(S_{2}) \subseteq (S_{2})_{II}$ and it follows that the non-trivial group element $(4,x^{3},2)(2,1,4)=(4,x^{3},4)$ belongs to $(S_{2})_{II}$. Therefore, $S_2$ does not have a one-point flow.}

\item{We show that $\RLM(S_{2})c=1$. Note that the 0-minimal ideal $I(S_{2})$ of $\RLM(S_{2})$ is the aperiodic 0-simple semigroup with structure matrix $M_{8}$. Therefore $\RLM(S_{2})c = (\RLM(S_{2})/I(S_{2}))c$ by the fundamental lemma of complexity \cite[Chapter 4]{qtheory}. Since every element in this Rees quotient acts by partial 1-1 maps on $I(S_{2})$, we have that $\RLM(S_{2})/I(S_{2})$ is a subsemigroup of the Symmetric Inverse Semigroup on $B$. Since every inverse semigroup has complexity at most 1, we are done, given that the element $a$ generates a non-trivial group.}

\end{enumerate}
 
In the following computation we use the notion of well-formed formulas, States and the evaluation transformation semigroup $Eval(S)$  of a $\GM$ semigroup $S$ from Section 5.1 of \cite{Trans}. We have summarized this material in the Appendix to the paper.

We now compute some states in the evaluation transformation semigroup $Eval(S_{2})$. We use the following notation for $SPC$. If
$(Y, \{\pi_{1}, \ldots , \pi_{n}\}, \{f_{1}, \ldots ,f_{n}\})$ is an $SPC$ we write $Y/\langle f_{1}|f_{2} \ldots |f_{n} \rangle$, where the separators
 correspond to the classes $\{\pi_{1}, \ldots \pi_{n}\}$. Thus in the example below, $1357/ \langle 1|1|1|1 \rangle $ is the $SPC$ with set $\{1,3,5,7\}$ and the identity partition with its unique cross-section. $2468/\langle 1xx^{2}x^{3} \rangle$ is the $SPC$ with set $\{2,4,6,8\}$, the partition with 1 block and cross section with representative $f$ with $f(2)=1,f(4)=x,f(6)=x^{2},f(8)=x^{3}$.  Consider the following sequence of computations that we let the reader justify.


$(1/\langle1\rangle)a^{\omega+*}=1357/\langle 1|1|1|1\rangle$, $(1357/\langle 1|1|1|1\rangle)b = 2/\langle 2\rangle$, $(2/\langle2\rangle)a^{\omega+*}=2468/\langle 1|1|1|1\rangle$, $(2468/\langle 1|1|1|1\rangle )b = 1\langle 1\rangle $.

It follows that if we let $D=(b(a^{\omega+*}))^{\omega+*}$, then $(1/\langle 1 \rangle )aD=12345678/\langle 11111111 \rangle$, where we have silently used Lemma 4.14.29 of \cite{qtheory} known as the Tie-Your-Shoes Lemma to show that we have one block. Therefore $\sigma_{4}=12345678/\langle 11111111\rangle$ is a State as are 
$\sigma_{1}=(12345678/\langle 11111111 \rangle)s_{x}=2468/\langle 1xx^{2}x^{3} \rangle,\sigma_{2}=(12345678/\langle 11111111\rangle)s_{x^{2}}=2468/\langle 1x^{2}1x^{2}\rangle$ and 
$\sigma_{3}=(12345678/\langle 11111111 \rangle )s_{x^{3}}=2468/ \langle 1x^{3}x^{2}x\rangle$. 

We now define a flow $F:(\{1,2,3,4\},RZ(4)^{1}) \rightarrow Rh_{B}(S_{2})$. We let $F(i)=\sigma_{i}, i =1,2,3,4$. We cover $a$ by the identity of $RZ(4)^{1}$, $b$ by $\bar{4}$, the constant map to $4$ and $s_{x^{i}}$ by $\bar{i}, i=1,2,3$. Finally we cover any element in $I(S_{2})$ by $\bar{4}$.

We verify that $F$ is a flow. We first claim that for $i=1,2,3,4$, $\sigma_{i}a=\sigma_{i}$. This is clear for $i=4$, since $a$ permutes $1 \ldots 8$ and acts as the identity on the $Z_4$ coordinate. We have $\sigma_{1}a=(2468/\langle 1xx^{2}x^{3}\rangle)a=2468/\langle x^{3}1xx^{2}\rangle$. But the cross-section $\langle x^{3}1xx^{2}\rangle =\langle 1xx^{2}x^{3}\rangle$ since $x(x^{3}1xx^{2})=(1xx^{2}x^{3})$. Therefore $F_{1}a = \sigma_{1}a=\sigma=F_{1}$. A similar calculation shows that $F_{2}a=F_{2}$ and $F_{3}a=F_{3}$. This is in fact using the elementary character theory of $Z_{4}$ as noted above. Therefore the conditions for a flow are verified for $a$. For $i=1,2,3$, 
$F(4)\sigma_{x^{i}}=((12345678/\langle 11111111\rangle)s_{x^{i}}=\sigma_{i}=F_{i}$ and is undefined on $F_{1},F_{2}$ and $F_{3}$. The image of of  $(a,g,b) \in I((S))$ on $F_{i}$ is either empty or less than or equal to $F_{4}$. This completes the proof that $F$ is a flow.

We therefore have that $(Z_{4} \times \{1,2,3,4\},S_{2}) \prec (Z_{4} \wr Sym(4) \wr (\{1,2,3,4\},RZ(4)^{1}) \times \RLM(S_{2})$ by Theorem \ref{PLFDT} and thus $S_{2}c=1$.
\ee

We now generalize the previous example $S_2$ to $S_n$ for all $n>0$.

\be

We modify the definition of $S_{2}$ in Example \ref{S2} as follows.

\begin{enumerate}
\item{Let $\Gamma_{2^{n+1}}$ be the incidence matrix of a $2^{n+1}$-cycle and let $M_{2^{n+1}}$ be the matrix $\begin{bmatrix} \Gamma_{2^{n+1}} & | & I_{2^{n+1}} \end{bmatrix}$ where $I_{2^{n+1}}$ is the identity matrix of size $2^{n+1}$.}

\item{Let $a$ be the permutation $(135 \ldots 2^{n+1}-1)(246\ldots 2^{n+1})$, so that $a$ is adding 2 modulo $2^{n+1}$. Let $b$ send $2^{n+1}$ to $1$ and $2^{n+1}-1$ to 2, so that $b$ sends an element of its domain from one orbit of $a$ to the other.}

\item{Let $C_{2^{n}}$ be the abstract character table of $Z_{2^{n}}$. With rows and columns labeled by $x^{j}, j = 0\ldots 2^{n}-1$, we have 
$C_{2^{n}}(x^{k},x^{l})=x^{kl}$, where we take $kl$ modulo $2^{n}$. Let $\chi_{i}$ be the character associated to $x^{i}$. Define $s_{i}$ to be the row monomial function sending $j$ to $\chi_{i}(j)2j$ for $j=1,\ldots 2^{n}$.}

\item{Define $S_{n}$ to be the subsemigroup of the translational hull of 
$I=M^{0}(Z_{2^n},\{1,\ldots 2^{n+2}\},\{1,\ldots 2^{n+1}\}, M_{2^{n+1}})$ generated by $a,b,s_{i}, i=1 \ldots 2^{n}$ and all elements of $I$. That these elements are in the translational hull of $I$ follows from \cite{cremona, deg2part2}.}

\item{Calculations analogous to those in Example \ref{S2} show that $S_{n}$ admits a flow over $RZ_{2^{n}}^1$ and that $S_{n}c=1$.}

\item{The resolution semigroup of the above flow decomposition is a small completely regular monoid having $Z_{2^{n}}$ as group of units  with generator the diagonal matrix with $(1,x,\ldots , x^{2^{n}-1})$ on the diagonal acting on the minimal ideal 
that is the simple semigroup with maximal group $Z_{2^{n}}$ and structure matrix the abstract character table $C_{2^n}$.}

\end{enumerate}
\ee
\section{Conclusion and Further Research}

In this paper we laid down the foundations of the theory of aperiodic flows for $\GM$ semigroups. In particular we gave a unified version of the  of the Presentation Lemma and its relation to flows and the Slice Theorem that have appeared in the literature. We characterized those $\GM$ semigroups that have flows over the trivial transformation semigroup as those that divide a direct product of a group and their $\RLM$ image. We showed how important results from inverse semigroup theory fit into the theory of 1-point flows. We looked at small monoids from the point of view of flows. We concluded with examples of semigroups of complexity 1 built from character tables of Abelian groups and built flows for them. In the follow up paper \cite{FlowsII} we prove new results in the theory of aperiodic flows and give more complicated examples. In \cite{FlowsIII} we develop the theory of non-aperiodic flows and give many examples. 

The document \cite{MasterList} contain many more examples of flows and examples in which aperiodic flows does not exist. The list will be updatged from time to time and the reader is advised to look at the list.

\appendix

\section{APPENDIX The Flow Monoid and the Evaluation Transformation Semigroup}

\subsection{The Monoid of Closure Operations}

In the Appendix we gather definitions and results from \cite{Trans}. For more details the reader is strongly urged to consult this paper. Since all the computations we do on examples in this paper are done in the Evaluation Transformation Semigroup defined in Section 5 of \cite{Trans} our goal is to summarize the background material in that paper needed to define this object. We begin with the definition of the monoid of closure operations on the direct product $L^{2} = L \times L$ of a lattice $L$ with itself. 

Let $L$ be a lattice and let $L^{2} = L \times L$. Let $f$ be a closure operator on $L^2$. By definition this means that $f$ is an order preserving, extensive (that is, for all $(l_{1},l_{2}) \in L^{2}, (l_{1},l_{2}) \leqslant (l_{1},l_{2})f$), idempotent function on $L^2$. A {\em stable pair} for $f$ is a closed element of $f$. Thus a stable pair $(l,l') \in L^2$ is an element such that $(l,l')f = (l,l')$.  The stable pairs of $f$ are a meet closed subset of $L \times L$. Conversely, each meet closed subset of $L \times L$ is the set of stable pairs for a unique closure operator on $L \times L$. We will identify $f$ as a binary relation on $L$ whose pairs are precisely the stable pairs. Let $B(L)$ be the monoid of binary relations on $L$. As is well known, $B(L)$ is isomorphic to the monoid $ M_{n}(\mathcal{B})$ of $n \times n$ matrices over the 2-element Boolean algebra $\mathcal{B}$, where $n=|L|$. The Boolean matrix associated to $f$ is of dimension $|L| \times |L|$. It has a 1 in position $(l_{1},l_{2})$ if $(l_{1},l_{2})$ is a stable pair and a 0 otherwise. 

With this identification, the collection $\mathcal{C}(L^{2})$ of all closure operators on $L^2$ is a submonoid of the monoid $B(L)$ of binary relations on $L$~\cite[Proposition 2.5]{Trans}. We thus also consider $\mathcal{C}(L^{2})$ to be a monoid of $|L| \times |L|$ Boolean matrices. Let $L$ be either $Rh_{B}(G)$ or $\operatorname{SP}(G\times B)$. We now recall the definition of some important unary operations on $\mathcal{C}(L^{2})$. 
\bd

Let Let $f \in \mathcal{C}(L^{2})$.

\begin{enumerate}

\item{The domain of $f$ denoted by $\operatorname{Dom}(f)$ is the set $\{x\mid \exists y, (x,y) \in f\}.$}

\item{Define the relation $\overleftarrow{f}$ by $\overleftarrow{f}= \{(x,x)\mid x \in \operatorname{Dom}(f)\}$. $\overleftarrow{f}$ is called {\em back flow along $f$}. See \cite[Remark 2.26]{Trans} for the reason for this terminology.}

\item{The relation $f^{*}$ is defined by $f^{*} = f \cap \{(x,x)\mid x \in X\}$. $f^{*}$ is called {\em the Kleene closure of $R$}. See Section 2 and Section 4 for the reason for this terminology.}

\item{Define the {\em loop of $f$} to be the relation $f^{\omega+*}=f^{\omega}f^{*}$, where $f^{\omega}$ is the unique idempotent in the subsemigroup generated by $f$.}

\end{enumerate}

\ed

At times it is convenient to identify $\overleftarrow{f}$ as $1|_{\operatorname{Dom}(f)}:L \rightarrow L$, the identity function restricted to 
$\operatorname{Dom}(f)$. Similarly, we identify $f^{*}$ as $1|_{\operatorname{Fix{f}}}:L \rightarrow L$, the restriction of the identity to the set of fixed-points of $f$, where $x \in L$ is a fixed point if $(x,x) \in f$. Our use of these will be clear from the context.

\subsection{The Flow Monoid}

Let $S$ be a $\GM$ semigroup generated by $X$ and let $x \in X$. Let $I(S) = M^{0}(G,A,B,C)$. We define a binary relation $f_{x}$
on $\operatorname{SP}(G \times B)$ by $((Y,\Pi), (Z,\Theta)) \in f_{x}$ if and only if $Yx \subseteq Z$ and the partial function induced by
right multiplication by $x$, $\cdot x: Y \rightarrow Z$ induces a well-defined partial injective map $\cdot x: Y/\Pi \rightarrow Z/\Theta$. This means that if  $(g,b), (g',b') \in Y$ and $(g,b)x,(g',b')x$ are both defined (and hence in $Z$ ), then $(g,b)\Pi (g',b')$ if and only if $(g,b)x\Theta (g',b')x$. Then $f_{x} \in \mathcal{C}(L^{2})$. See \cite[Proposition 2.22]{Trans}. $f_{x}$ is called the free-flow by $x$.

We now define the flow monoid $M(L)$ as follows. This is called the 0-flow monoid in \cite{Trans}, but since this paper only considers aperiodic flow, we just say flow monoid.

\bd

Let $L = \operatorname{SP}(G \times B)$. The flow monoid $M(L)$, is the
smallest subset of $\mathcal{C}(L^{2})$ satisfying the following axioms:

\begin{enumerate}

\item {(Identity) The multiplicative identity $I$ of $\mathcal{C}(L^{2})$ is in $M(L)$.}

\item{(Points) For all $x \in X$, $f_{x}$ the free-flow along $x$ belongs to $M(L)$.}

\item{(Products) If $f_{1}, f_{2} \in M(L)$, then $f_{1}f_{2} \in M(L)$.}

\item{(Vacuum) If $f \in M(L)$, then $\overleftarrow{f} \in M(L)$.}

\item{(Loops) If $f \in M(L)$, then $f^{\omega+*}\in M(L)$.}

\end{enumerate}

\ed

We remark that for each $n \geq 0$, there is an $n$-flow monoid $M_{n}(L)$ defined in \cite{Trans}. The definition above of $M(L)$ is denoted by $M_{0}(L)$ in \cite{Trans} and used extensively in \cite{complexity1}. For $n>0$, Axioms (1)-(4) are the same for $M_{n}(L)$ as for $M_{0}(L)$. Axiom (5) restricts the use of the loop operator to $n$-loopable elements \cite[Section 4]{Trans}. In \cite{complexityn}, $n$-loopable elements are replaced by a more restrictive definition and play a crucial role in the main results of \cite{complexityn}.

\subsection{The Evaluation Transformation Semigroup}

We now defined the Evaluation Transformation Semigroup $\mathcal{E}(L) = (\operatorname{States}, \operatorname{Eval}(L))$ as in \cite{Trans}. Again, because of our interest in this paper on aperiodic flows, we don't define the $n$-Evaluation Transformation Semigroups $E_{n}$ for $n>0$. We begin with the definition of Well-Formed Formualae (WFFs).

\bd

Let $X$ be an alphabet. We define a well-formed formula inductively as follows.

\begin{enumerate}

\item{The empty string $\epsilon$ is a well-formed formula.}

\item{Each letter $x \in X$ is a well-formed formula.}

\item{If $\tau, \sigma$ are well-formed formulae, then
so is $\tau\sigma$.}

\item{If $\tau$ is a well-formed formula that is not a proper power (i.e., not of the form
$\sigma^{n}$ where $n > 1$), then $\tau^{\omega+*}$ is also a well-formed formula.} 

\end{enumerate}

\ed

The set of well-formed
formulae is denoted by $\Omega(X)$. Well-formed formulae will be denoted by Greek
letters. As a convention, if $\tau=\sigma
^{n}$, where $\sigma$ is not a proper power, then we set
$\tau^{\omega+*}=\sigma^{\omega+*}$. In other words, we extract roots before applying the unary operation $\omega+*$.

Let $V = \prod_{f \in M(L)}\overleftarrow{f}$. $V$ is called the {\em Vacuum}. See \cite{Trans} for an explanation of this terminology. In \cite{Trans}, $V$ is denoted by $\mathcal{F}_{0}$. It is proved in \cite{Trans} that $V$ is an idempotent in $M(L)$. We now will work exclusively in the subsemigroup $VM(L)V$ of $M(L)$. We want to interpret WFFs in $VM(L)V$.

\bd

Define recursively a partial function $\mathcal{I}: \Omega(X)\rightarrow VM(L)V$ as follows. 

\begin{enumerate}

\item{$\epsilon\mathcal{I} = V$.}

\item{$x\mathcal{I} = VxV$ for $x \in X$.}

\item{If $\mathcal{I}$ is already defined on $\tau, \sigma \in \Omega(X)$, set $(\tau\sigma)\mathcal{I} = \tau\mathcal{I}\sigma\mathcal{I}$.}

\item{If $\tau \in \Omega(X)$ is not a proper power and $\tau\mathcal{I}$ is defined, set $\tau^{\omega+*}\mathcal{I} = 
(\tau\mathcal{I})^{\omega+*}$.}

\end{enumerate}

\ed

We normally omit $\mathcal{I}$ and assume that a WFF $\tau$  is being evaluated in $VM(L)V$ according to the definition of $\mathcal{I}$. 
 We first define a new operator on elements of $M(L)$ called {\em forward flow}. Recall that the bottom of the lattice $\operatorname{SP}(G \times B)$ is the pair $\Box = (\emptyset, \emptyset)$.

\bd

Let $f \in M(L)$ and let $l \in L$. Let $(l,\Box)f =(l_{1},l_{2})$ We define the forward flow of $f$ denoted by $\overrightarrow{f}$ by 
$l\overrightarrow{f}=l_{2})$. That is, we apply $f$ to $(l,\Box)$ and project to the right-hand coordinate.

\ed

The following is proved in Sections 2 and 4 of \cite{Trans}.

\begin{enumerate}

\item{$\overrightarrow{f}:L \rightarrow L$ is an order preserving function on $L$.}

\item{If $lV=l$, then $(l,\Box)f=(l,l')$ and $l' \in LV$. That is, if $l \in LV$ and $f \in VM(L)V$ then 
$l\overrightarrow{f} \in LV$. Furthermore, the left-coordinate of $(l,\Box)f$ is still $l$ and the right-hand coordinate is in $LV$. Therefore $\overrightarrow{f}$ is a well-defined function from $LV$ to $LV$.}

\item{The assignment of $f$ to $\overrightarrow{f}$ defines an action of $VM(L)V$ on $LV$. It follows that we have a transformation semigroup $(LV,M'(L))$, where $M'(L)$ is the image of $VM(L)V$ on $LV$ under this action.}

\end{enumerate}

We now restrict the action in $(LV,M'(L))$ to the set of $\operatorname{States}$ defined as follows. Let $(g,b) \in G \times B$. Then 
the element $(\{(g,b)\},\{(g,b)\}) \in L$ is called a {\em point}. In Section 5 of \cite{Trans}, it is proved that every point $p$ satisfies $pV=p$. 

\bd

The set of $\operatorname{States}$ is the smallest subset of $LV$ such that:

\begin{enumerate}

\item{Every point $p \in \operatorname{States}$.}

\item{If $l \in \operatorname{States}$, then $l\overrightarrow{f} \in \operatorname{States}$.}

\end{enumerate}

\ed

In other words $\operatorname{States}$ is the smallest subset of $LV$ containing the points and closed under the action of $M(L)$ on $LV$. We can finally define the Evaluation Transformation Semigroup where all the computations in the examples in this paper take place.

\bd

The Evluation Transformation Semigroup $\mathcal{E}(L)$ is defined by $\mathcal{E}(L) = (\operatorname{States}, \operatorname{Eval}(L))$ where $\operatorname{Eval}(L)$ is the image of $M'(L)$ by restricting its action to $\operatorname{States}$.

\ed

We remark that there is an Evaluation Transformation Semigroup $\mathcal{E}_{n}(L)$ for all $n \geq 0$. The definition here is the case $n=0$
which is all we need in this paper as we are only concerned with aperiodic flows.

\bibliography{stubib}
\bibliographystyle{abbrv}


\end{document}